\documentclass[11pt]{article}
\usepackage{xcolor}
\usepackage{amsmath}
\usepackage{amssymb}
\usepackage{amsthm}
\usepackage{latexsym}
\usepackage[all,cmtip]{xy}
\usepackage{graphicx, float, color}
\usepackage{dsfont}
\usepackage[affil-it]{authblk}

\newcommand{\E}{\mathcal{E}}
\newcommand{\B}{\mathcal{B}}
\newcommand{\inv}{^{-1}}
\newcommand{\ch}{\mathds{1}}
\newcommand{\cspan}{\overline{\text{span}}}
\newcommand{\N}{\mathbb{N}}
\newcommand{\Z}{\mathbb{Z}}

\DeclareMathOperator{\ind}{ind}
\DeclareMathOperator{\res}{res}
\DeclareMathOperator{\Id}{Id}
\DeclareMathOperator{\spn}{span}
\DeclareMathOperator{\supp}{supp}
\DeclareMathOperator{\Ob}{Obj}

\theoremstyle{plain}
\newtheorem{thm}{Theorem}[section]
\newtheorem{definition}[thm]{Definition}
\newtheorem{lemma}[thm]{Lemma}
\newtheorem{proposition}[thm]{Proposition}
\newtheorem{corollary}[thm]{Corollary}

\newtheorem{lem}[thm]{Lemma}
\newtheorem{defn}[thm]{Definition}
\newtheorem{cor}[thm]{Corollary}

\theoremstyle{definition}
\newtheorem{example}[thm]{Example}
\newtheorem{rmk}[thm]{Remark}

\begin{document}

\title{Discrete Conduch\'{e} Fibrations and C*-algebras}

\author{Jonathan H. Brown}
\affil{Mathematics Department\\
University of Dayton\\
300 College Park\\
Dayton, OH 45469-2316\\
USA.}

\author{David N. Yetter}
\affil{Mathematics Department\\
    Kansas State University\\
138 Cardwell Hall\\
Manhattan, KS 66506-2602\\
USA.}

\maketitle

\begin{abstract}
The $k$-graphs in the sense of Kumjian and Pask \cite{KP} are discrete Conduch\'{e} fibrations over the monoid ${\mathbb N}^k$ in which ever morphism in the base has a finite preimage under the fibration.  We examine the generalization of this construction to discrete Conduch\'{e} fibrations with the same finiteness condition and a lifting property for completions of cospans to commutative squares, over any category satisfying a strong version of the right Ore condition, including all categories with pullbacks and right Ore categories in which all morphisms are monic.
\end{abstract}

\section{Introduction}
In 2000 Kumjian and Pask introduce $k$-graphs in order to generalize the high-rank Cuntz-Krieger algebras of Robertson and Steger \cite{RS} and the graph algebras of Kumjian, Pask, Raeburn and Renault \cite{KPRR}.   Their construction of $C^*$-algebras from $k$-graphs yield algebras which are uncommonly tractable and include a wide variety of examples (see for example \cite{PRRS}).  This combination has lead to significant interest in $k$-graph $C^*$-algebras.

A $k$-graph consists of a category $\Lambda$ that is fibred over $\N^k$ by a degree functor $d:\Lambda\to \N^k$ satisfying the unique factorization condition: if $d(\lambda)=m+n$ then there exist unique $\mu,\nu$ with $d(\mu
)=m, d(\nu)=n$ and $\lambda=\mu\nu$ \cite{KP}.  Kumjian and Pask suggest that  it might be interesting to study categories $\Lambda$ that are fibred over a cancellative abelian monoid.  Recently \cite{CKSS} and \cite{Yang} use such a construction as a technical tool to study the primitive ideals and abelian subalgebras of $k$-graphs respectively.  A category fibred over a cancellative abelian monoid (in particular a $k$-graph) is an example of a previously studied notion in category theory called a discrete Conduch\'e fibration (see Definition~\ref{def:dCf} below).  Our goal in this paper is to develop a theory of $C^*$-algebras associated to discrete Conduch\'e fibrations and thus generalize $k$-graph $C^*$-algebras.  While extensions of the theory of \cite{KP} to categories fibred over a cancellative abelian monoid is relatively straightforward, substantial difficulties arise when the base of the fibration is a more general category.

More recently, Spielberg \cite{S} has proposed a different construction of $C^*$-algebras from what he calls categories of paths.  These are special categories which, among other things, contain no inverses to non-identity morphisms and are both left and right cancellative.   He also does not require any fibration of the category corresponding to the degree functor and so his construction has a different flavor from \cite{KP}. We allow our categories to have inverses.   Thus  a discrete group $H$ equiped with the identity functor $\Id_H: H\to H$ is an admissible fibration for our theory but neither category is a category of paths in the sense of Speilberg.  It turns out that the $C^*$-algebra associated by our construction to  $\Id_H: H\to H$ is the full group $C^*$-algebra of $H$.   We also consider categories that aren't necessarily cancellative.  

In Section 2 we introduce a discrete Conduch\'e fibration $F:\E\to \B$  and show that there is a universal $C^*$-algebra, $C^*(F)$ for the analogous Cuntz-Kreiger relations.  In Section 3 we restrict our attention to fibrations over strongly right Ore categories (see Definition~\ref{def:strongly right ore}).  This restriction allows us to define an infinite path as a section of the functor induced by $F$ on a slice category.  Our infinite paths  behave analogously to those of \cite{KP} (see Proposition~\ref{restrictandinduce}).   Then in Proposition~\ref{prop:nonzero} we represent $C^*(F)$ on the set of infinite paths and use this to show $C^*(F)\neq \{0\}$.   In Section 4 we put a topology on the set of infinite paths and then use a special class of local homomorphisms to construct a locally compact Hausdorff \'etale groupoid $G_F$.  In Section 5 we show that if morphisms in the category are monic then $C^*(G_F)\cong C^*(F)$.  We use this to show that if morphisms are also epi then the category $\cal E$ can be represented faithfully as a subcategory of operators on a Hilbert space.

Throughout our discussion, all categories are small, and composition in categories is written in the customary anti-diagrammatic order as multiplication, thus the composition of the morphisms in $C \stackrel{\beta}{\rightarrow} B \stackrel{\alpha}{\rightarrow} A$ is written $\alpha\beta$.  We denote the target or codomain of a morphism $\alpha$ by $r(\alpha)$ and its source or domain by $s(\alpha)$.  Given an object $X$ in a category $\mathcal{C}$ we denote by $X\mathcal{C}$ the set of morphisms in $\mathcal{C}$ with range $X$.  We denote the set of objects in $\mathcal{C}$ by $\Ob(\mathcal{C})$, similarly for a morphism $\gamma$ in $\mathcal{C}$ we denote $\gamma\mathcal{C}=\{\gamma\mu:s(\gamma)=r(\mu),\mu\in \mathcal{C}\}$.  We will often identify categories with their set of morphisms, so that $\gamma \in {\mathcal {C}}$ will mean that $\gamma$ is an morphism of $\mathcal{C}$.

\section{Discrete Conduch\'e fibrations}

\begin{definition}
\label{def:dCf}
A {\em discrete Conduch\'{e} fibration (dCf)} is a functor $F:{\cal E}\rightarrow {\cal B}$ with the {\em unique factorization lifting property} or {\em Conduch\'{e} condition}:  for every morphism $\phi:Y\rightarrow X$ in ${\cal E}$, every factorization of $F(\phi)$ in ${\cal B}$

\[ 
F(Y)\stackrel{\lambda} {\rightarrow}B \stackrel{\rho}{\rightarrow} F(X)
 \]

\noindent  lifts uniquely to a factorization of $\phi$

\[ Y\stackrel{\tilde{\lambda}}{\rightarrow} Z \stackrel{\tilde{\rho}}{\rightarrow} X \]

\noindent with $F(\tilde{\lambda}) = \lambda$, $F(\tilde{\rho}) = \rho$, and $F(Z) = B$.
\end{definition}

As an aside we note that the more categorically natural notion of a Conduch\'{e} fibration, in which the factorization is unique only up to isomorphism appears to be inadequate for our purposes, though in most of our examples, there are no non-identity isomorphisms in $\cal E$, and in such cases the two notions coincide.

A simple argument shows that unique factorization extends to finitely many factors:

\begin{lemma}  \label{morefactors} If $F:{\cal E}\rightarrow {\cal B}$ is a discrete Conduch\'{e} fibration, $\mu:Y\rightarrow X$ is any morphism of ${\cal E}$ and $F(\mu)=abc$, that is

\[ F(Y) \stackrel{a}{\rightarrow} B_1 \stackrel{b}{\rightarrow}  B_2 \stackrel{c}{\rightarrow}  F(X) \]

\noindent is a factorization of $F(\mu):F(Y)\rightarrow F(X)$ in ${\cal B}$, then there exists a unique
factorization

\[ Y \stackrel{\alpha}{\rightarrow} Z_1 \stackrel{\beta}{\rightarrow}  Z_2 \stackrel{\gamma}{\rightarrow}  X \]

\noindent of $\mu$ in ${\cal E}$ such that $F(\alpha ) = a$, $F( \beta) = b$, $F( \gamma) = c$,
$F(Z_1) = B_1$, and $F(Z_2) = B_2$.

Moreover, in general, any factorization of a morphism in ${\cal B}$ into a finite number of composands lifts uniquely to a factorization of any morphism in its preimage in ${\cal E}$ into the same number of composands.
\end{lemma}

\begin{lemma} \label{identitieslift} If $F:{\cal E}\rightarrow {\cal B}$ is a discrete Conduch\'{e} fibration, and $F(\phi:X\rightarrow Y) = \Id_B$ for some object $B$ in $\cal B$, then $\phi$ is an identity morphism in $\cal E$.
\end{lemma}

\noindent {\sc Proof.} By functoriality of $F$ we have $F(X) = F(Y) = B$.  Now $\Id_B$ factors as $\Id_B(\Id_B)$, and  both $\Id_Y(\phi)$ and $\phi(\Id_X)$ are lifts of this factorization.  Therefore, by the uniqueness condition of dCFs, $\Id_Y = \phi = \Id_X$, and we see $\phi$ is an identity morphism. $\Box$ \medskip

Following Kumjian, Pask and Raeburn \cite{KPR} we define:

\begin{definition}
A functor $F:{\cal E}\rightarrow {\cal B}$ is {\em row-finite} if for every object $X$ in ${\cal E}$ and every morphism $\beta:B\rightarrow F(X)$ in $\cal B$, the class of morphisms with target
$X$ whose image under $F$ is $\beta$ is a finite set.
\end{definition}

The following surjectivity condition takes the place of the no sources condition in the study of graphs and $k$-graphs.

\begin{definition} A functor $F:{\cal E}\rightarrow {\cal B}$ between small categories is {\em strongly surjective} if it is surjective on objects, and given any object $X$ in ${\cal E}$ the map induced from the set of morphisms targetted at $X$  to the set of morphisms targetted at $F(X)$ in ${\cal B}$ is surjective.
\end{definition}

We may then rephrase the definition of \cite{KP} as

\begin{definition}
A {\em $k$-graph} is a countable category ${\cal E}$ equipped with a strongly surjective, row-finite, dCf to the additive monoid ${\mathbb N}^k$, regarded as a category with one object.
\end{definition}

Now it is easy enough to describe generators and relations after the manner of Cuntz-Krieger \cite{CK}, Kumjian, Pask and Raeburn \cite{KPR}, and Kumjian and Pask \cite{KP} using an arbitrary strongly surjective, row-finite functor between small categories as data:

\begin{definition} \label{CKsystem}
Given a strongly surjective, row-finite functor $F:{\cal E}\rightarrow {\cal B}$, a {\em Cuntz-Krieger system associated to $F$} is the following data in a $C^*$-algebra $D$:

\begin{itemize}
\item A projection $P_X$ for each object $X$ of ${\cal E}$,
\item A partial isometry $S_\alpha$ for each morphism $\alpha:Y\rightarrow X$ of ${\cal E}$,
\end{itemize}

\noindent satisfying the relations
\begin{enumerate}
\item For $X \neq Y$,  $P_X \perp P_Y$,
\item If $\alpha$ and $\beta$ are composable, then $S_{\alpha\beta} = S_\alpha S_\beta$,
\item For all $X$, $P_X = S_{\Id_X} = S^*_{\Id_X}$,
\item For all $\alpha:Y\rightarrow X$, $S^*_\alpha(S_\alpha) = P_Y$,
\item If $F(\alpha) = F(\beta)$ and $\alpha \neq \beta$ then $S^*_\beta(S_\alpha) = 0$,
\item For all $X$, and for all morphisms $b:B\rightarrow F(X)$ in ${\cal B}$
\[ \sum_{\{\alpha \in X{\cal E} : F(\alpha) = b\}} S_\alpha S^*_\alpha = P_X. \]
\end{enumerate}

We denote the smallest $C^*$-subalgebra of $D$ containing $P,S$ by $C^*(P,S).$
\end{definition}

We show in the next proposition that there exists a universal $C^*$-algebra for these relations. The proof of this proposition follows the argument of \cite[Theorem~2.1]{HR}.

\begin{proposition}
\label{c*alg}
Let $F:{\cal E}\to {\cal B}$ be a discrete Conduch\'e Fibration. Then there exists a $C^*$-algebra, $C^*(F)$ generated by a Cuntz-Krieger system $\{p_X, s_\alpha\}$ such that for any Cuntz-Krieger system $\{Q_{X}, T_\alpha\}$, associated to $F$ in a $C^*$-algebra $B$ there is a unique $*$-homomorphism from $C^*(F)$ to $B$ extending the map $s_\alpha\mapsto T_\alpha$.
\end{proposition}

\noindent {\sc Proof.} 
 Let $K_F$ be the free complex  $*$-algebra generated by the morphisms of ${\cal  E}$: that is $K_F$ is a vector space over $\mathbb{C}$ with basis given by the set of words ${\mathcal{W}}$ in symbols $\{\alpha:\alpha\in \E\}\cup\{\beta^*:\beta\in \E\}\cup \{X: X\in \Ob(\E)\}$, multiplication induced by concatenation of words and involution induced by the maps $\alpha\in \E \mapsto \alpha^*$ and $c\in \mathbb{C}\mapsto\overline{c}$. Define a norm $\|\cdot \|_I$ on $K_F$ by $\|\sum_{w\in N} c_w w\|_I=\sum_{w\in N}|c_w|$ where $N$ is a finite subset of ${\cal W}$ and $c_w$ are complex numbers.

Let $J$ be the ideal in $K_F$ generated by the Cuntz-Krieger relations.  Then $\Gamma(F)=K_F/J$ and $\|\cdot\|_I$ induces a norm $\|\cdot\|_\Gamma$ on $\Gamma(F)$ by 
\[
\|a+J\|_\Gamma=\inf_{b+J=a+J} \|b+J\|_I.
\]

\noindent Let 
\[
\|a\|_0=\sup\{ \|\rho(a)\| : \rho \text{~is a $\Gamma$-norm decreasing $*$-representation of $\Gamma(F)$ on a Hilbert space}\}
\]
for $a\in \Gamma(F)$.  Since $\|\cdot\|_0$ is bounded by $\|\cdot\|_\Gamma$ this gives a $C^*$-seminorm on $\Gamma(F)$.

Now define $C^*(F)_0:=\Gamma(F)/\ker (\|\cdot\|_0)$.  Then $\|\cdot\|_0$ induces a $C^*$-norm $\|\cdot\|$ on $C^*(F)_0$.  Complete $C^*(F)_0$ in this norm to obtain a $C^*$-algebra $C^*(F)$. Denote the images of objects $X$ and morphisms $\mu$ of {\cal E} in $C^*(F)$ by $p_X$ and $s_\mu$ respectively.   

It remains to show that this $C^*$-algebra satisfies the required universal property.   Suppose $\{Q_X\}_{X\in \Ob(\mathcal{E})}$ and $\{T_\mu\}_{\mu\in \mathcal{E}}$ is a Cuntz-Krieger family in a $C^*$-algebra $B$.  By the Gelfand-Nainmark theorem we can assume $B$ is a subalgebra of the bounded operators on a Hilbert space.  By the universal property of the free algebra $K_F$ there exists a $*$-homomorphism $\tilde{\rho}$ from $\Gamma(F)$ to $B$.  By the definition of $\|\cdot\|_0$ this representation is norm decreasing on $C^*(F)_0$ and thus extends to $*$-homomorphism $\rho$ from $C^*(F)$ to $B$ as desired.  By construction $\rho$ must send $p_X$ to $Q_X$ and $s_\mu$ to $T_\mu$ and thus is determined on the dense subalgebra $C^*(F)_0$ and thus $\rho$ is unique. $\Box$
\smallskip

A question immediately arises:  Under what circumstances is the map $s_{(\cdot )}:{\cal E}\rightarrow
C^*(F)$ injective?

It is the upshot of results in \cite{KP} is that when $F$ is a row-finite, strongly surjective dCF over ${\mathbb N}^k$, i.e. a $k$-graph, the map $s_{(\cdot )}:{\cal E}\rightarrow
C^*(F)$ is injective.  This result uses an infinite path construction.  We will use a similar construction below to prove Proposition~\ref{prop:nonzero} which gives some mild conditions that ensure that $C^*(F)\neq \{0\}$ and  Proposition~\ref{prop:inj} which gives conditions that guarantee that the map ${\cal E}\to C^*(F)$ given by $\alpha\mapsto s_\alpha$ is injective.

\section{ Kumjian-Pask Fibrations and Infinite Paths}

Recall that a morphism $\alpha$ is  {\em epi} if $\beta\alpha=\gamma\alpha$ implies $\beta=\gamma$; $\alpha$ is {\em monic} if $\alpha\beta=\alpha\gamma$ implies $\beta=\gamma$. A category $\cal B$ is {\em right (resp. left) cancellative} if every morphism is epi (resp. monic).

One initially might think that cancellation conditions after the manner of Spielberg \cite{S} are necessary. However, as we shall see, weaker conditions suffice for most of the construction, though for some results cancellation conditions on the base category are required.

Following Johnstone \cite{J} define:

\begin{definition}
\label{def:strongly right ore}
A category {\em satisfies the right Ore condition} or for brevity is {\em a right Ore category} if every cospan $A\stackrel{m}{\rightarrow} B \stackrel{n}{\leftarrow} C$
can be completed to a commutative square.
\end{definition}

\noindent We extend this notion to

\begin{definition} \label{strongOre} A category is {\em strongly right Ore} if it is right Ore and moreover, for every cospan $m,n \in \B$ with $p_1, p_2, q_1, q_2\in \B$ with $mp_i=nq_i$ for $i=1,2$ there exist $a,b\in \B$ with $p_1a=p_2b$ and $q_1a=q_2b$.

\end{definition}

Such categories are plentiful:

\begin{proposition} 

~ 

\begin{enumerate}
\item Any category with pullbacks is strongly right Ore.
\item Any left cancellative right Ore category is strongly right Ore.
\end{enumerate}
\end{proposition}

\noindent{\sc Proof.} For (1),  the existence of pullbacks implies that the category is right Ore.  To see such a category is strongly right Ore, suppose $(P,t,w)$ is the pullback of the cospan $(m,n)$ so that $mt=nw$ and $(Z_1, p_1, q_1), (Z_2,p_2, q_2)$ are two competions of $(m,n)$ so that $mp_i=nq_i$.  By the universal property of pullbacks there exists $u_i: Z_i\to P$ such that $tu_i=p_i$ and $wu_i=q_i$.  Now $(u_1,u_2)$ is a cospan so let $(Q,a,b)$ be the pullback of this cospan.  By definition we have $p_1a=p_2b$ and $q_1a=q_2b$. 

For (2), given a cospan $m,n$, let $mp_i=nq_i$ for $i=1,2$ be two completions to commutative squares, and let $a$ and $b$ complete the cospan $mp_1, mp_2$ to a commutative square
$mp_1a = mp_2b$.  Now, since $m$ is monic, it follows that $p_1a = p_2b$.  Since the cospan completed could also be considered to be $nq_1, nq_2$, and $n$ is monic, it follows that $q_1a = q_2b$. $\Box$

\smallskip

Since  all lattices in the order-theoretic sense (including the lattices of open sets in topological spaces) have pull-backs they are all strongly right Ore.  Similarly since all groups, groupoids and the positive cones of lattice-ordered groups are left cancellative they are strongly right Ore as well.

\begin{proposition}
If ${\cal C}_i$ for $i\in {\cal I}$ is a set-indexed family of categories, each of which is strongly right Ore, then the product category, $\prod_{i \in {\cal I}} {\cal C}_i$, is strongly right Ore.
\end{proposition}

\noindent{\sc Proof.}  The required completions of diagrams exist in the product because their components exist in the factors.  $\Box$
\smallskip

\begin{definition} \label{slice}

Given a category ${\cal C}$ and an object $X$ in ${\cal C}$, the {\em slice category} ${\cal C}/ X$ is the category whose objects are morphisms of ${\cal C}$ with target $X$, and morphisms given by commutative triangles:  that is for $\alpha,\beta\in X{\cal C}=\text{Obj}({\cal C})/X$, a morphism from $\beta$ to $\alpha$ is given by a $\gamma\in {\cal C}$ with $\alpha\gamma=\beta$.   We can thus view morphisms in ${\cal C}/ X$ as ordered pairs $(\alpha,\gamma)\in X{\cal C}\times {\cal C}$  with $r(\gamma)=s(\alpha)$.  The range of $(\alpha,\gamma)$ is $\alpha$ and its source is $\alpha\gamma$.  Composition $(\alpha,\gamma)(\beta,\delta)$ can only occur if $\beta=\alpha\gamma$ and in this case we have $(\alpha,\gamma)(\alpha\gamma,\delta)=(\alpha, \gamma\delta)$.  Let $\pi_i$ be the projection onto the $i$th factor of ${\cal C}/ X$. Notice that $\pi_1(\alpha,\gamma)\pi_2(\alpha,\gamma)=\alpha\gamma$ and  that $\pi_2((\alpha,\gamma)(\alpha\gamma,\delta))=\pi_2(\alpha,\gamma)\pi_2(\alpha\gamma,\delta)$. (Slice categories are a special instance of a more general construction known as a comma category. See \cite{CWM} for more details. )

 \end{definition}

For any object $X$ in ${\cal C}$,  a functor $F:{\cal C}\rightarrow {\cal D}$ will induce a natural functor  from ${\cal C}/ X$ to ${\cal D}/ F(X)$ (which by abuse of notation we also denote $F$).

\begin{definition}
A functor $F:{\cal E}\rightarrow {\cal B}$ is {\em locally split} if for every object $X \in {\cal E}$, the induced functor (also denoted $F$ by abuse of notation) $F:{\cal E}/ X \rightarrow {\cal B}/ F(X)$ admits a splitting (or section) $x:{\cal B}/ F(X) \rightarrow {\cal E}/ X$, that is a functor such that $F\circ x = \Id_{{\cal B}/ F(X)}$. (Contrary to the usual practice in category theory of denoting functors by capital letters, we will use lower case Latin letters near the end of the alphabet to denote local splittings of the fibration, as these will play the role of infinite paths, denoted by such letters in the work of Kumjian and Pask.)
\end{definition}

Note that a locally split functor which is surjective on objects is strongly surjective. 

Let $x:{\cal B}/ F(X) \rightarrow {\cal E}/ X$ be a splitting of $F$ as above and $(a,b)\in {\cal B}/ F(X)$. Then $a$ is an object in ${\cal B}/ F(X)$. By unique factorization,  for $x(a,b)=(x_1(a,b),x_2(a,b))$,   we have $x(a)=x_1(a,b)=x_2(\Id_X, a)$.  

\begin{definition} 
A functor $F:{\cal E}\rightarrow {\cal B}$ targetted at a strongly right Ore category is a {\em 
Kumjian-Pask fibration (KPf)} if it is a locally split dCf.
\end{definition}

Notice that we do not require strong surjectivity, however we have

\begin{proposition}
If $F:{\cal E}\rightarrow {\cal B}$ is a KPf, then so is $\hat{F}:{\cal E}\rightarrow F({\cal E})$. Moreover
$\hat{F}$ is strongly surjective.
\end{proposition}

\noindent {\sc Proof.} 
Since $F$ is locally split dCF, $\hat{F}$ is too.  $\hat{F}$ is surjective on objects and therefore strongly surjective.  To show $\hat{F}$  is a KPf it suffices to show $F({\cal E})$ is strongly right Ore.   We begin by showing that $F({\cal E})$ is right Ore.  Let $(m,n)$ be a cospan in $F({\cal E})$ and $A=r(m)$.  Since $m\in  F({\cal E})$ we can choose $\mu\in {\cal E}$ with $F(\mu)=m$.  Thus $F(r(\mu))=r(m)=A$.  Let $X=r(\mu)$.  Since $F$ is locally split, choose a section $x: {\cal B}/ A\to {\cal E}/ X$.  Since $m,n\in {\cal B}$ and ${\cal B}$ is right Ore there exists $p,q\in {\cal B}$ such that $mp=nq$.  Now $F(x_2(m,p))=p$ and $F(x_2(n,q))=q$ and thus $p,q\in F({\cal E})$: that is $F({\cal E})$ is right Ore.  Similarly if $mp_i=nq_i$ then since ${\cal B}$ is strongly right Ore, there exist $a,b\in {\cal B}$ with $p_1a=p_2b$ and $q_1a=q_2b$.  So $F(x_2(mp_1,a))=a$ and $F(x_2(np_2,b))=b$ and so $a,b\in F({\cal E})$, that is $F({\cal E})$ is strongly right Ore as well. $\Box$

There are many examples:  

\begin{example}  If $d:\Lambda \rightarrow {\mathbb N}^k$ is a $k$-graph, then $d$ is a KPf.
\end{example}

\begin{example} If $\cal B$ is any strongly right Ore category, $\Id_{\cal B}:{\cal B}\rightarrow {\cal B}$ is a KPf.
\end{example}

Of particular interest is a special case of the last example:

\begin{example} If $H$ is a discrete group, regarded as a one object category, then the identity group homomorphism, $\Id_H$, regarded as a functor is a KPf.  In this example the $C^*$-algebra $C^*(\Id_H)$ is particularly easy to compute. Suppose $(S,P)$ is a Cuntz-Krieger system for $\Id_H$ in a $C^*$-algebra $A$. Since $H$ has only one object $e$ and we have $P_e S_t=S_tP_e$ for all $t\in H$, and so $(S,P)$ is a Cuntz-Krieger system in the unital $C^*$-algebra $PAP$ with unit $P$.  Because we are using the identity functor the sum in relation (6) has only one summand and, combining with relation (4), we have $S_tS_t^*=P=S_t^*S_t$ for all $t\in H$.  This gives that $S_t$ is unitary.  Now by relation (2) we get that $t\mapsto S_t$ is a group homomorphism, and thus is a unitary representation of $H$.  Since $C^*(H)$ is universal for unitary representations of $H$ and $C^*(\Id_H)$ is universal for Cuntz-Krieger systems we get that $C^*(\Id_H)\cong C^*(H)$.
\end{example}

As a preliminary to our next examples, recall that  any poset $(P, \leq)$ induces a category, which by abuse of notation we also denote $P$, whose objects are elements of $P$ and whose morphisms are ordered pairs $(p,q)\in P\times P$ with $p\leq q$, with $r((p,q)) = q$ and $s((p,q)) = p$ and composition defined by $(q,r)(p,q)=(p,r)$ (cf. \cite{CWM}).  Notice that since reversing the ordering gives another poset, there is an arbitary choice in this convention, which was made so that the  abstract maps of a poset would correspond in the case of concrete posets of subobjects of a mathematical object to the inclusion maps from one subobject to another which contains it.
(The convention, unfortunately, appears to be at odds with that adopted by Kumjian and Pask \cite{KP} for their category $\Omega_k$ in which the morphisms have the greater $k$-tuple as source and the lesser as target.  In fact, $\Omega_k$ is the opposite category of $({\mathbb N}^k, \leq)$ regarded as a category in the standard way, a fact either intentionally or fortuitously emphasized by the fact Kumjian and Pask denote the unique morphism from $p^\prime$ to $p$ when $p^\prime \geq p$ by $(p,p^\prime)$, the same notation by which the standard convention for a poset would denote the unique morphism in the opposite direction.  As we will see in Example \ref{k-graphs}, $\Omega_k$ is better thought of as arising not from the poset of $k$-tuples of natural numbers, but as the slice category of the monoid $({\mathbb N}^k, +, 0)$ (as a one-object category) over its unique object, which naturally gives the morphisms in the direction from greater to lesser, though it will also suggest use of a different notation for them.)

\begin{example} Let $X$ be a topological space, and $\mathfrak{S}$ a presheaf of sets on $X$; that is a contravariant functor from the lattice of open sets $\mathcal{ T}_X$ ordered by inclusion to ${\bf Sets}$, the category of sets and set-functions.  As usual, we regard the images of, elements of $\mathfrak{S}(U)$ under maps $\mathfrak{S}(V,U)$ as restrictions, denoting $\mathfrak{S}(V,U)(\sigma)$ by $\sigma |_V$ for $\sigma \in \mathfrak{S}(U)$.

Let $\mathcal{S}$ be the set of local sections, that is pairs $(U,\sigma)$ where $U$ is an open set and $\sigma \in {\mathfrak S}(U)$
Then $\mathcal{S}$ is partially ordered by restriction, with $(V,\tau) \leq (U,\sigma)$ exactly when $V\subseteq U$ and $\tau = \sigma |_V$.  Now $((V,\tau),(U,\sigma))\mapsto (V,U)$ defines a functor $F: {\cal S}\to {\cal T}_X$.  If $((V,\tau), (U,\sigma))\in {\cal S}$ and $V\subseteq W\subseteq U$ so that $(V,W)(W,U)$ is a factorization of $F((V,\tau), (U,\sigma))=(V,U)$, then $((V,\tau), (W,\sigma |_W))((W,\sigma |_W), (U,\sigma))$ is the unique factorization of $((V,\tau), (U,\sigma))$ lifting the factorization of $(V,U)$ as $(V,W)(W,U)$.   

Since ${\cal T}_X$ has pullbacks it is strongly right Ore.  And, finally, for any object $(U, \sigma)$ in $\mathcal{S}$, the induced functor $F: \mathcal{S}/(U,\sigma) \rightarrow \mathcal{T}_X/U$ admits a splitting given by $(V,U) \mapsto ((V,\sigma |_V),(U,\sigma))$, and thus $F:\mathcal{S}\rightarrow \mathcal{T}_X$ is a KPf.
\end{example}

In fact, this last example generalizes.

\begin{example}  If $\cal B$ is a right Ore poset (i.e. a poset in which every pair of elements which admits an upper bound admits a lower bound, necessarily strongly right Ore since all morphisms are monic) and $S:{\cal B}^{op} \rightarrow {\bf Sets}$ is a presheaf of  sets on $\cal B$, then the poset of sections ${\cal E}_S$ with elements all pairs $(U, a)$ with $U \in Ob({\cal B})$ and $a \in S(U)$, and $(V, b) \leq (U, a)$ exactly when $V \leq U$ and $b = a|_V$ (using the usual notation in the context of (pre)sheaves on a space) admits an obvious monotone map $(U, a) \mapsto U$ to $\cal B$, which, regarded as a functor, is a KPf.
\end{example}

Among the last two classes of examples, those arising from presheaves of finite sets will be row-finite.

\begin{definition}
\label{infinite paths}
Let $F:{\cal E}\rightarrow {\cal B}$ be a row-finite, Kumjian-Pask fibration. 
\begin{enumerate}
\item For an object $X$ of ${\cal E}$, {\em an infinite path to $X$} is a section $x: {\cal B}/ F(X) \rightarrow {\cal E}/ X$ of the functor $F:{\cal E}/ X \rightarrow {\cal B}/ F(X)$. 
\item  We denote the set of all infinite paths to $X$ by $Z(X)$, and the set of all infinite paths (that is $\bigcup_{X\in Ob({\cal E})} Z(X)$) by $F^\infty$.
\item For a morphism $\mu:Y\rightarrow X$ in ${\cal E}$, {\em an infinite path ending in $\mu$} is an infinite path $x$ to $X$ for which $x(F(\mu)) = \mu$.  We denote the set of all infinite paths ending in $\mu$ by $Z(\mu)$.
\item When an infinite path $x$ is introduced without explicitly stating that it is a functor from ${\cal B}/F(X)$ to ${\cal E}/X$ for a specified object $X$, we will denote the object over which its target category is a slice category by $r(x)$.
\end{enumerate}
\end{definition}

These notions are consistent with the usual notion of infinite paths from row-finite source free $k$-graphs.  Indeed the following example is the motivation for regarding splittings of the induced functors on slice categories as infinite paths.

\begin{example} \label{k-graphs} Let $d:\Lambda \rightarrow {\mathbb N}^k$ be a $k$-graph (for some $k > 0$).  The base category is the additive monoid ${\mathbb N}^k$, regarded as a category with a single object $\star$ in the usual way.

Using the conventions of Definition \ref{slice}, the objects of  ${\mathbb N}^k/\star$ are elements of ${\mathbb N}^k$, while morphisms are ordered pairs $(p,q)\in {\mathbb N}^k \times {\mathbb N}^k$, with range $p$ and source $p+q$.  Note that if $p$ and $p^\prime$ are two objects of ${\mathbb N}^k/\star$ there is a unique
morphism from $p^\prime$ to $p$ exactly when $p^\prime \geq p$ in the product partial order on ${\mathbb N}^k$, and thus, the slice category is isomorphic to  the category Kumjian and Pask \cite{KP} denote by $\Omega_k$.  The unique map from $p^\prime$ to $p$, when one exists, which Kumjian and Pask denote by $(p,p^\prime)$ is denoted in our convention for slice categories by $(p, p^\prime - p)$.

The functor by which Kumjian and Pask endow $\Omega_k$ with a $k$-graph structure maps all objects to $\star$ and the morphism which we denote by $(p, p^\prime - p)$ and Kumjian and Pask denote $(p, p^\prime)$ to $p^\prime - p$.  This is precisely the natural forgetful functor from the slice category ${\mathbb N}^k/\star$ to
${\mathbb N}^k$.

It remains to show that an infinite path in the sense of Kumjian and Pask \cite{KP}, that is a degree preserving functor from $\Omega_k$ to $\Lambda$, or equivalently from ${\mathbb N}/\star$ to $\Lambda$, is equivalent to a splitting of the induced functor $d:\Lambda/v \rightarrow {\mathbb N}/\star$.  Given any functor $F:{\cal C}\rightarrow {\cal D}$ and an object $X$ of $\cal C$ the square formed by the given $F$, the induced $F:{\cal C}/X \rightarrow {\cal D}/F(X)$ and the forgetful functors from the slice categories to the ${\cal C}$ and ${\cal D}$ commutes.  

Applying this observation to the degree functor of $\Lambda$, we see that if $\tilde{x}$ is a splitting of
$d:\Lambda/v \rightarrow  {\mathbb N}/\star$, and $U:\Lambda/v \rightarrow \Lambda$ is the forgetful functor, then $x:=U\tilde{x}$ is a degree-preserving functor from ${\mathbb N}/\star$ to $\Lambda$, since 
$d U \tilde{x} = \Upsilon d \tilde{x} = \Upsilon$, where $\Upsilon:{\mathbb N}/\star \rightarrow {\mathbb N}$ is the forgetful functor, which coincides with Kumjian and Pask's degree functor on $\Omega_k$ by the preceding discussion, the first equation holding by the commutativity of the square of functors just observed, the second by the fact that $\tilde{x}$ splits the induced $d$.  In particular $x$ will be an infinite path in the sense of \cite{KP} to the object $v$ of $\Lambda$.

Conversely, if $x:{\mathbb N}^k\rightarrow \Lambda$ is a degree preserving functor, we construct a factorization of $x$ of the form $U\tilde{x}$ for $\tilde{x}:{\mathbb N}^k/\star \rightarrow \Lambda/x(0)$ a splitting
of the induced $d:\Lambda/x(0) \rightarrow {\mathbb N}^k/\star$ as follows:

For an object $q$ of ${\mathbb N}^k/\star$, that is a morphism of ${\mathbb N}^k$, let $\tilde{x}(q) = x(q,0)$.
For a morphism $(p,q)$ of ${\mathbb N}^k/\star$ from $p+q$ to $q$, let $\tilde{x}(p,q) = (x(p,q),x(q,0))$.  The functoriality of $\tilde{x}$ follows easily from that of $x$.  It follows from the preservation of degree by $x$ that
$d \tilde{x}(q) = q$ and $d \tilde{x}(p,q) = d(x(p,q), x(q,0)) = (p,q)$, so that $\tilde{x}$ is a splitting of the induced
functor $d$ on slice categories.   Finally since the forgetful functor from a slice category to the underlying category is given in our notation by second projection $U \tilde{x}(q) = U x(q, 0) = x(q)$ and $U \tilde{x}(p,q) = U(x(p,q), x(q,0)) = x(p,q)$.

It is easy to see that the two constructions are inverse to each other.

\end{example}

The following is easy to verify

\begin{lemma} \label{partitionbylifts}

If $b:B\rightarrow F(X)$ is a morphism in ${\cal B}$, then

\[ Z(X) = \bigcup_{\{ \beta \in{X\cal E} :  \text{~and~} F(\beta) = b\}} Z(\beta) ,\]

\noindent and the union is disjoint.
\end{lemma}

The hypothesis that $F$ is  locally split does not seem to guarantee there be any infinite paths ending in a particular morphism. However, the following proposition in particular implies that for a Kumjian-Pask fibration $F:\E\to \B$ and $\alpha\in \E$, there exists an infinite path $x\in Z(\alpha)$.  This proposition defines maps which are crucial to the rest of our study so we enumerate their properties here.

\begin{proposition} \label{restrictandinduce}
If $F:{\cal E}\rightarrow {\cal B}$ is a Kumjian-Pask fibration, for any morphism $\mu:Y\rightarrow X \in {\cal E}$, there exist bijections
\begin{align*}
\ind_\mu &:Z(Y) \rightarrow Z(\mu), \quad\text{and}\\
\res_\mu &:Z(\mu) \rightarrow Z(Y)
\end{align*}
satisfying the following properties:
\begin{enumerate}
\item\label{it: ind further} $\ind_\mu(x)(F(\mu)a)=\mu x(a)$ for any $a\in F(s(\mu))\B$;
\item\label{it: inverse} $\ind_\mu \circ \res_\mu=\Id_{Z(\mu)}$ and $\res_\mu\circ \ind_\mu=\Id_{Z(r(\mu))}$;
\item\label{it: ind res on objects} for $X$ an object in $\E$, $\ind_{\Id_X}=\res_{\Id_X}=\Id_{Z(X)}$;
\item\label{it:resres} for $\mu,\nu\in \E$ with $s(\mu)=r(\nu)$, $\res_\mu(Z(\mu\nu))= Z(\nu)$ and $\ind_\mu(Z(\nu))=Z(\mu\nu)$;
\item\label{it:res homo} for $\mu,\nu\in \E$ with $s(\mu)=r(\nu)$, then the domain of $\res_\nu\circ\res_\mu$ is $Z(\mu\nu)$ and $\res_\nu\circ\res_\mu=\res_{\mu\nu}$;
\item\label{it:ind homo} for $\mu,\nu\in \E$ with $s(\mu)=r(\nu)$, then the domain of $\ind_\mu\circ\ind_\nu$ is $Z(s(\nu))$ and $\ind_\mu\circ\ind_\nu=\ind_{\mu\nu}$.
\end{enumerate}
\end{proposition}

\noindent {\sc Proof.} We begin by defining the maps $\ind_\mu$ and $\res_\mu$.  Since $\res_\mu$ is easier to describe we define it first.  

Let $x$ be an infinite path ending in $\mu$, $x:{\cal B}/ F(Y) \rightarrow {\cal E}/ Y$.  Recall that for $(c,d)\in  {\cal B}/ F(Y)$, we denote $x(c,d)=(x_1(c,d),x_2(c,d))$ and then $x(cd)$ is a morphism in $\E$ whose factorization with respect to $cd$ is $x_1(c,d)x_2(c,d)$.  The idea of $\res_\mu(x)$ is that it removes the first $\mu$ part of the morphism $x(F(\mu)ab)$.

As objects can be identified with identity morphisms, it is enough to describe $\res_\mu(x)$  on pairs of the form $(a,b)$ with $r(a)=F(s(\mu))$ and $r(b)=s(a)$.  Now $(F(\mu)a,b)\in {\cal B}/ F(Y) $ and
\begin{align*}
x(F(\mu) a,b)&=(x_1(F(\mu) a,b),x_2(F(\mu) a,b))=(x(F(\mu)a), x_2(F(\mu) a,b))\\
&=(\mu x_2(F(\mu),a), x_2(F(\mu) a,b))\in \E/ Y\subset {\cal E}\times {\cal E}.
\end{align*}  
We then define
\[
\res_\mu(x)(a,b):=(x_2(F(\mu),a),x_2(F(\mu) a,b)).
\]
As $r(x_2(F(\mu),a))=s(x(F(\mu)))=s(\mu)$ and $r(x_2(F(\mu) a,b))=s(x(F(\mu)a))=s(x_2(F(\mu),a))$ we get that $\res_\mu(x)$ maps into ${\cal E}/ s(\mu)$.  As $x$ is a section,
 \[
F(\res_\mu(x)(a,b))=(F(x_2(F(\mu),a)),F(x_2(F(\mu) a,b)))=(a,b),\] so $\res_\mu(x)$ is a section of $F: {\cal E}/s(\mu)\to {\cal B}/F(s(\mu))$.  Lastly, to see $x$ is a functor,  if $(a,b),(ab,c)\in {\cal B}/ F(s(\beta))$, then 

\begin{align*}
\res_\mu(x)(a,b)\res_\mu(ab,c)&=(x_2(F(\mu),a),x_2(F(\mu) a,b))(x_2(F(\mu),ab),x_2(F(\mu) ab,c))\\
&=(x_2(F(\mu),a), x_2(F(\mu) a,b)x_2(F(\mu) ab,c)))\\
&=(x_2(F(\mu),a), x_2(F(\mu) a,b)(x(F(\mu) ab,c)))\\
&=(x_2(F(\mu),a), x_2(F(\mu) a,bc))=\res_\mu(x)(a,bc).
\end{align*}
Thus $\res_\mu(x):\B/ F(s(\mu))\to \E/ s(\mu)$ as desired.

We turn our attention to defining $\ind_\mu$.  The idea of $\ind_\mu$ is to add $\mu$ onto the begining of infinite paths. Suppose $x\in Z(s(\mu))$.  Consider $(a,b)\in {\cal B}/ F(Y)$.  Now $F(\mu), ab$ is a cospan in ${\cal B}$.  Since ${\cal B}$ is right Ore there exists $(c,d)\in {\cal B}$ such that $F(\mu)c=abd$.  By unique factorization there exist unique morphisms $\alpha^{c,d}_{a,b},\beta^{c,d}_{a,b}, \delta^{c,d}_{a,b}\in {\cal E}$ such that 
\[
\mu x(c)=\alpha^{c,d}_{a,b}\beta^{c,d}_{a,b}\delta^{c,d}_{a,b}
\]
and $F(\alpha^{c,d}_{a,b})=a, F(\beta^{c,d}_{a,b}(b))=b$, and $F(\delta^{c,d}_{a,b}(d))=d$. Take 
\[
\ind_\mu(x)(a,b)=(\alpha^{c,d}_{a,b},\beta^{c,d}_{a,b}).
\]   

We need to show this does not depend on the choice of the span $c,d$.  For this we use the strongly right Ore condition.  Suppose there exist $c_i,d_i$ that complete the cospan $F(\mu), ab$ into a commuting square.  By strongly right Ore, there exist $t,u\in {\cal B}$ such that $c_1t=c_2 u$ and $d_1  t=d_2 u$.  Thus it suffices to show that $(\alpha^{c_1,d_1}_{a,b},\beta^{c_1,d_1}_{a,b})=(\alpha^{c_1t,d_1t}_{a,b},\beta^{c_1t,d_1t}_{a,b})$.  Consider $\mu x(c_1 t)=\mu x(c_1)x_2(c_1, t)$.   By definition 
\[
\mu x(c_1)x_2(c_1, t)=\mu x(c_1t)=\alpha^{c_1t,d_1t}_{a,b}\beta^{c_1t,d_1t}_{a,b}\delta^{c_1t,d_1t}_{a,b}
\]
with $F(\alpha^{c_1t,d_1t}_{a,b})=a, F(\beta^{c_1t,d_1t}_{a,b})=b$ and $F(\delta^{c_1t,d_1t}_{a,b})=d_1t$. By unique factorization there exist $\gamma, \tau$ with $\gamma\tau=\delta^{c_1t,d_1t}_{a,b}$ and $F(\gamma)=d_1$ and $F(\tau)=t$.  Since $F(\alpha^{c_1t,d_1t}_{a,b}\beta^{c_1t,d_1t}_{a,b}\gamma)=abd_1=F(\mu x(c_1))$, by unique factorization we have
\[
\alpha^{c_1,d_1}_{a,b}\beta^{c_1,d_1}_{a,b}\delta^{c_1,d_1}_{a,b}=\mu x(c_1)=\alpha^{c_1t,d_1t}_{a,b}\beta^{c_1t,d_1t}_{a,b}\gamma
\]
so by unique factorization again we have $\alpha^{c_1,d_1}_{a,b}=\alpha^{c_1t,d_1t}_{a,b}$ and $\beta^{c_1,d_1}_{a,b}=\beta^{c_1t,d_1t}_{a,b}$  as desired.  Thus $\ind_\alpha(x)$ is well-defined.  

We need to check that $\ind_\mu(x)$ is a section of $F:{\cal E}/ Y\to {\cal B}/ F(Y)$.  Notice that $r(\alpha_{a,b}^{c,d})=r(\mu)$ and $r(\beta_{a,b}^{c,d})=s(\alpha_{a,b}^{c,d})$ by definition so $\ind_\mu(x)(a,b)\in {\cal E}/ Y$.  Also by definition $F(\ind_\mu(x)(a,b))=(a,b)$.  It remains to show that $\ind_\mu(x)$ is a functor.  Suppose $(a,b)(ab,t)\in  {\cal B}/ F(Y)$.  Consider the cospan $F(\mu)$ and $abt$.  By right Ore, there exist $c,d\in {\cal B}$ with $F(\mu)c=abtd$.  So by unique factorization there exist unique $\alpha,\beta,\tau,\delta$ such that $F(\alpha)=a, F(\beta)=b, F(\tau)=t$ and $F(\delta)=d$ and $\mu x(c)=\alpha\beta\tau\delta$.  By definition then

\begin{align*}
\alpha^{c,td}_{a,b}&=\alpha=\alpha^{c,d}_{a,bt},\\
\beta^{c, td}_{a,b}&=\beta,\\
\alpha^{c,d}_{ab,t} &=\alpha\beta,\\
\beta^{c,d}_{ab,t}&=\tau,\\
\beta^{c,d}_{a,bt}&=\beta\tau\\
\end{align*}
Thus
\[
\ind_\mu(x)(a,b)\ind_\mu(x)(ab,t)=(\alpha, \beta)(\alpha\beta,\tau)=(\alpha,\beta\tau)=\ind_\mu(a,bt) \]
so that $\ind_\mu$ is indeed a functor.


To see that $\ind_\mu(x)\in Z(\mu)$ and that $\res_\mu$ is the inverse for $\ind_\mu$ we need the following claim: 
\begin{equation}
\label{clm1}
 (\ind_\mu(x))_1 (F(\mu)a,b) = \mu x(a)\quad\text{and}\quad (\ind_\mu(x))_2 (F(\mu) a, b) = x_2(a,b).
 \end{equation} 
To see this note that $F(\mu), F(\mu)ab$ is a cospan and the span $ab, s(b)$ completes it to a commuting square.  So we have by the definition of $\ind_\mu(x)$ that 
\[
\mu x(a)x_2(a,b)=\mu x(ab) = (\ind_\mu(x))_1(F(\mu)a,b) (\ind_\mu(x))_2(F(\mu)a,b).
\]
Unique factorization now gives equation~\eqref{clm1}.

Since $\ind_\mu(x)(F(\mu)a) = \ind_\mu(x)_1(F(\mu)a,F(s(\mu)))=\mu x(a)$,  item~\eqref{it: ind further} follows immediately from equation~\eqref{clm1}; in particular, taking $a=s(\alpha)$ gives $\ind_\mu(x)\in Z(\mu)$. 

We now show item~\eqref{it: inverse}, that is $\ind_\mu$ is the inverse of $\res_\mu$.  Let $x\in Z(\mu)$ and $(a,b)\in {\cal B}/ F(r(\mu))$; we want to show $\ind_\mu(\res_\mu(x))(a,b)=x(a,b)$.  We compute $\ind_\mu(\res_\mu(x))(a,b)$.  Now $F(\mu), ab$ is a cospan in ${\cal B}$ so there exist $c,d\in {\cal B}$ with $F(\mu)c=abd$. 
Thus
\begin{align*}
x_1(a,b)x_2(a,b)x_2(ab,d) &=x_1(ab,d)x_2(ab,d)\\
&=x(abd)\\
&=x(F(\mu)c)\\
&=\mu x_2(F(\mu),c)\quad\text{since }x\in Z(\mu)\\
&=\mu\res_\mu(x)(c)\quad\text{by definition of }\res_\mu\\
&=(\ind_\mu(\res_\mu(x)))_1(a,b) (\ind_\mu(\res_\mu(x)))_2(a,b)\delta_{a,b}^{c,d}
\end{align*}
by definition of $\ind_\mu$.  Thus $\ind_\mu(\res_\mu(x))(a,b)=x(a,b)$ by unique factorization. 

Now suppose that $x\in Z(s(\mu))$  and $(a,b)\in {\cal B}/ F(s(\mu))$; we want to show $\res_\alpha (\ind_\mu(x))(a,b)=x(a,b)$. 
Using equation~\eqref{clm1} we compute:
 \begin{align*}
\res_\mu(\ind_\mu(x)(a,b))&=(\ind_\mu(x)_2(F(\mu),a), \ind_\mu(x)_2(F(\mu)a,b)))\\
&=(x_2(s(\mu),a),x_2(a,b))=(x_1(a,b),x_2(a,b))=x(a,b)
\end{align*}
as desired. 

The other properties now follow fairly quickly.

For item~\eqref{it: ind res on objects}, note that 
\[
\res_{\Id_X}(x)(a,b)=(x_2(\Id_X,a),x_2(\Id_X a,b))=(x_1(a,b),x_2(a,b))=x(a,b).
\]
So that $\res_{\Id_X}=\Id_{Z(X)}$.  Now $\ind_{\Id_X}=(\res_{\Id_X})\inv=\Id_{Z(X)}\inv=\Id_{Z(X)}$.

For item~\eqref{it:resres}, we first show 
\begin{equation}
\label{inclusion}
\ind_\mu(Z(\nu))\subset Z(\mu\nu)\quad\text{and}\quad\res_\mu(Z(\mu\nu))\subset Z(\nu).
\end{equation}
For the former note if $x\in Z(\nu)$, $\ind_\mu(x)(F(\mu)F(\nu))=\mu x(F(\nu))=\mu\nu$ by item~\eqref{it: ind further}. For the later, if $x\in Z(\mu\nu)$, then $\mu\nu=x(F(\mu)F(\nu))=x_1(F(\mu),F(\nu))x_2(F(\mu),F(\nu))$. Since
$F(x_2(F(\mu),F(\nu)))=F(\nu)$ we have $x_2(F(\mu),F(\nu))=\nu$ by unique factorization.  Therefore 
\[ \res_{\mu}(x)(F(\nu))= \res_\mu(x)_1 (F(\nu),s(\nu)))=x_2(F(\mu),F(\nu))=\nu ;\]

\noindent that is $\res_{\mu}(x)\in Z(\nu)$.

It now suffices to show the reverse inclusions.  If $x\in Z(\nu)$, then $\ind_\mu(x)\in Z(\mu\nu)$.  So $x=\res_\mu(\ind_\mu(x))\in \res_\mu(Z(\mu\nu))$; that is $Z(\nu)\subset  \res_\mu(Z(\mu\nu))$ and so $\res_\mu(Z(\mu\nu))=Z(\nu)$.  Similarly, if $x\in Z(\mu\nu)$, then $\res_\mu(x)\in Z(\nu)$, so that $x=\ind_\mu(\res_\mu(x))\in \ind_\mu(Z(\nu))$; that is $Z(\mu\nu)\subset \ind_\mu(Z(\nu))$ and so $\ind_\mu(Z(\nu))=Z(\mu\nu)$.

For item~\eqref{it:res homo}, we first show that the domain of $\res_\nu\circ\res_\mu=Z(\mu\nu)$.  Let $x$ be  in the domain of $\res_\nu\circ\res_\mu$, we what to show $x\in Z(\mu\nu)$. Now  $x=\ind_\mu\circ\ind_\nu(  \res_\nu\circ\res_\mu(x))$ so that 
\[
x(F(\mu)F(\nu))=\ind_\mu\circ\ind_\nu(  \res_\nu\circ\res_\mu(x))(F(\mu)F(\nu))=\mu \ind_\nu(  \res_\nu\circ\res_\mu(x))(F(\nu))=\mu\nu
\]
where we used item~\eqref{it: ind further} twice.  Thus $x\in Z(\mu\nu)$.  Now if $x\in Z(\mu\nu)$ then by item~\eqref{it:resres}, $\res_\mu(x)\in Z(\nu)$ and so $x$ is in the domain of $\res_\nu\circ \res_\mu$ as desired.  

Next we show $\res_\nu\circ\res_\mu=\res_{\mu\nu}$ are equal.  For this we compute
\begin{align*}
\res_\nu\circ &\res_\mu(x)(a,b)=((\res_\mu(x))_2 (F(\nu),a), (\res_\mu(x))_2(F(\nu)a,b))\\
&=(x_2(F(\mu\nu),a),x_2(F(\mu\nu)a,b))\\
&=\res_{\mu\nu}(x)(a,b)
\end{align*}
as desired.

For item~\eqref{it:ind homo}, notice that by item~\eqref{it:res homo} we have $\res_\nu\circ\res_\mu=\res_{\mu\nu}$, so 
\[
\ind_\mu\circ\ind_\nu=(\res_\nu\circ\res_\mu)\inv=(\res_{\mu\nu})\inv=\ind_{\mu\nu}
\]
finishing the proof.
$\Box$ \medskip

\begin{rmk}
Notice that to define $\res_\beta$ only requires that $F:{\cal E}\to {\cal B}$ is a dCF.  We only use that $F$ is a KPf to define $\ind_\beta$.
\end{rmk}

We can now show that $C^*(F)\neq 0$ for Kumjian-Pask fibrations.

\begin{proposition}
\label{prop:nonzero}
Let $F:\E\to \B$ be a row-finite Kumjian-Pask fibration.  Then $C^*(F)\neq 0$.
\end{proposition}

\noindent {\sc Proof.}  Consider $\ell^2(F^\infty)$.  For $\mu\in \E$ define 
\[
T_\mu x=\begin{cases} \ind_\mu x & \text{if~} x\in Z(s(\mu))\\
 0 &\text{otherwise}
 \end{cases}
 \]
 and $Q_X$ to be the projection onto the subspace spanned by $Z(X)$.  A quick computation shows that 
 \[
T_\mu^* x=\begin{cases} \res_\mu x & \text{if~} x\in Z(\mu)\\
 0 &\text{otherwise}.
 \end{cases}
 \]
 
Proposition~\ref{restrictandinduce} shows that $Q,T$ satisfies all of the conditions of a Cuntz-Krieger $F$-family except (6).  But (6) follows imediately from Lemma~\ref{partitionbylifts}.  Thus, by the universal property,  $C^*(Q,T)$ is a quotient of $C^*(F)$.  Since $C^*(Q,T)$ is nonzero we have that $C^*(F)\neq 0$.
$\Box$
\section{The groupoid of a Kumjian-Pask fibration}

In this section we use the infinite path space of a Kumjian-Pask fibration to construct a groupoid.

A groupoid $G$ is a small category in which every morphism is invertible.  We identify the objects in $G$ with the identity morphisms and denote both by $G^{(0)}$.  As $G$ is a category we can send any morphism $\gamma$ to its range and source and denote these maps by $r$ and $s$ respectively.  A topological groupoid is a groupoid with a topology in which both composition is continuous and inversion is a homeomorphism.  It follows that $r$ and $s$ are also continuous. An open subset $B$ of a topological groupoid $G$ is called a bisection if $r|_B$ and $s|_B$ are injective and open. We say a topological groupoid is \'etale if it has a basis of bisections.  

We are interested in locally compact Hausdorff \'etale groupoids because there is a well developed theory of the $C^*$-algebras constructed from them.  We will show that Kumjian-Pask fibrations give rise to \'etale groupoids,  but first we sketch the construction of locally compact Hausdorff \'etale groupoid $C^*$-algebras for the convenience of the reader (for details see \cite{ren80}).  Let $G$ be a locally compact Hausdorff \'etale groupoid.  Define a convolution algebra structure on the continuous compactly supported functions $C_c(G)$ on $G$ by
\[ f*g(\gamma)=\sum_{\eta : r(\eta)=r(\gamma)} f(\eta)g(\eta\inv \gamma)\quad \quad \quad f^*(\gamma)=\overline{f(\gamma\inv)}:\]
the sum is finite because $f$ is compactly supported and the inverse image of a point under an \'{e}tale map (here $r^{-1}(\gamma)$) is discrete. We say a sequence of functions $\{f_i\}$ in $C_c(G)$ converges to $f$ in the inductive limit topology if $\|f_i-f\|_\infty=\sup_{\gamma\in G}\{|f_i(\gamma)-f(\gamma)|\}\to 0$ and there exists a compact set $K\subset G$ such that $\supp(f_i)\subset K$ eventually.   Define  $\text{Rep}(C_c(G))$ to be the set of $*$-homomorphisms from $C_c(G)$ to $B(\mathcal{H})$ such that the image of inductive limit convergent sequences are weak-$*$ convergent.  Define a norm on $C_c(G)$ by
\begin{align*}
\|f\| &=\sup\{\|\pi(f)\|: \pi\in \text{Rep}(C_c(G))\}
\end{align*}
We then define $C^*(G)$ to be the completion of $C_c(G)$ in $\|\cdot\|$.(It is non trivial that the  supremum defining the norm exists, see \cite{ren80}.) Given a unit $u\in G^{(0)}$ there is a representation $L^u: C^*(G)\to B(\ell^2(Gu))$ given by 

\[
L^u(f)\delta_{\xi}=\sum_{s(\eta)=r(\xi)} f(\eta)\delta_{\eta\xi}.
\]
A quick computation shows that $L^u\in \text{Rep}(C_c(G))$.  We define $I_\lambda=\cap_{u\in G^{(0)}} \ker(L^u)$ and $C^*_r(G)=C^*(G)/I_\lambda$.

To describe the groupoid constructed from a Kumjian-Pask fibration we first need to topologize the infinite path space. In Definition~\ref{infinite paths} we described for each $\alpha\in {\cal E}$ a set of infinite paths $Z(\alpha)=\{x\in F^\infty: x(F(\alpha))=\alpha\}$.  In this section we show that, under a mild countability hypothesis, the collection of these sets forms a basis of compact open sets for a locally compact Hausdorff topology on $F^\infty$.

\begin{lem}
\label{inclusion exclusion}
Let $F:{\cal E}\rightarrow {\cal B}$ be a row-finite Kumjian-Pask fibration
\begin{enumerate}
\item If $\alpha,\delta\in {\cal E}$ with $r(\delta)=s(\alpha)$ then $Z(\alpha\delta)\subset Z(\alpha)$.
\item If $\alpha,\beta\in {\cal E}$ with $F(\alpha)=F(\beta)$ then $Z(\alpha)\cap Z(\beta)\neq \emptyset$ if and only if $\alpha=\beta$.
\end{enumerate}
\end{lem}

\noindent {\sc Proof.} For the first item observe that if $x\in Z(\alpha\delta)$ then $x(F(\alpha)F(\delta))=\alpha\delta$ so by unique factorization $x(F(\alpha))=\alpha$.

For the second item if $x\in Z(\alpha)\cap Z(\beta)$ then $\alpha=x(F(\alpha))=x(F(\beta))=\beta$.
$\Box$ \medskip

\begin{proposition}
\label{disjoint union}
Let $F:{\cal E}\rightarrow {\cal B}$ be a  Kumjian-Pask fibration and $\alpha,\beta\in {\cal E}$ with $Z(\alpha)\cap Z(\beta)\neq \emptyset$.  Then there exists $c,d\in {\cal B}$ such that $F(\alpha)c=F(\beta)d$, 
$I=\{\alpha\gamma: s(\alpha)=r(\gamma), F(\gamma)=c\}\cap \{\beta\delta: s(\beta)=r(\delta), F(\delta)=d\}\neq\emptyset$ and 
\[
Z(\alpha)\cap Z(\beta)=\bigcup_{\mu\in I} Z(\mu).
\]
Further the sets $Z(\mu)$ for $\mu\in I$ are mutually disjoint.
\end{proposition}

\noindent {\sc Proof.} If $x\in Z(\alpha)\cap Z(\beta)$ then $r(\alpha)=r(x)=r(\beta)$.  Thus $F(\alpha), F(\beta)$ is a cospan in ${\cal B}$ and since ${\cal B}$ is right Ore, there exist $c,d\in {\cal B}$ with $F(\alpha)c=F(\beta)d$.  For $x\in Z(\alpha)\cap Z(\beta)$,  $\alpha x_2(F(\alpha),c)=x(F(\alpha)c)=x(F(\beta)d)=\beta x_2(F(\beta),d)$, so we have $\alpha x_2(F(\alpha),c)\in I$  and $x\in Z(\alpha x_2(F(\alpha),c))$. Thus $Z(\alpha)\cap Z(\beta)\subset \bigcup_{\mu\in I} Z(\mu)$.  Now the reverse inclusion follows from Lemma~\ref{inclusion exclusion}.  Finally  Lemma~\ref{inclusion exclusion} also gives that  the $Z(\mu)$ are mutually disjoint. 
$\Box$ \medskip

\begin{corollary}
\label{basis for F}
Let $F:{\cal E}\rightarrow {\cal B}$ be a  Kumjian-Pask fibration.  Then the set $\{Z(\alpha):\alpha\in {\cal E}\}$  is a basis for a  topology on $F^\infty$.  Furthermore, if $\cal B$ is countable and $F$ is row-finite, then
$F^\infty$ is second countable.  If each slice category ${\cal B}/ B$ is countable and $F$ is row-finite then each $Z(\beta)$ and each $Z(X)$ is second countable.
\end{corollary}

\noindent {\sc Proof.} Since $F^\infty=\bigcup_{X\in \Ob({\cal E})} Z(X)$, Proposition~\ref{disjoint union} shows $\{Z(\alpha):\alpha\in {\cal E}\}$  is a basis.   For the second statement, if $\B$ is countable and $F:\E\to \B$ is row-finite, then $\E$ is countable and  $\{Z(\alpha):\alpha\in {\cal E}\}$ is countable. $\Box$ \medskip

\begin{rmk}
$\B$ countable implies that ${\cal B}/ B$ is countable for all objects $B$ in $\B$.  The later can happen without $\B$ being countable for instance if $\B$ is the uncountable disjoint union of countable categories.  
\end{rmk}

Henceforth we regard $F^\infty$ as a topological space with the topology induced by the basis $\{Z(\beta) : \beta \in {\cal E}\}$, and the subsets $Z(\beta)$ and $Z(X)$ for morphisms $\beta$ and object $X$ of $\cal E$ as spaces in the subspace topology.

The following lemma will be quite useful in what follows:

\begin{lemma} \label{objectssuffice}
Let $F:\E\to \B$ be a KPf, $X$ an object in $\E$, $x,y\in Z(X)$ and $(a,b)\in \B/ F(X)$. The following are equivalent

\begin{enumerate}
\item $x(a,b) = y(a,b)$  as morphisms in ${\cal E}/ F$.
\item $x(a) = y(a)$ and $x(ab) = y(ab)$  as objects in ${\cal E}/ F$.
\item $x(ab) =y(ab)$ as objects in ${\cal E}/ F$.
\end{enumerate}

\noindent Consequently if $x$ and $y$ agree on all objects of $B/ F(X)$ then $x = y$.
\end{lemma}

\noindent{\sc Proof.}   For the first statement, note that (1) implies (2) since equality of morphisms implies equality of their sources and targets, and (2) implies (3)
trivially.  To see that (3) implies (1), recall that $x_1(a,b)x_2(a,b)=x(ab)=y(ab)=y_1(a,b)y_2(a,b)$.  Thus by unique factorization $x_i(a,b) = y_i(a,b)$ for $i = 1,2$ and therefore $x(a,b)=y(a,b)$. 

The second statement follows from the first, since agreement on the source of a morphism implies agreement on the morphism. $\Box$
\medskip

\begin{proposition}
Let $F:\E\to \B$ be a KPf.  Then $F^\infty$ is Hausdorff.
\end{proposition}


\noindent{\sc Proof.}  
Suppose $x\in Z(X)$ and $y\in Z(Y)$ with $x\neq y$.  If $X\neq Y$ then $Z(X)\cap Z(Y)=\emptyset$ and we are done.  If $X = Y$, then by Lemma \ref{objectssuffice} there is a morphism (object in the slice category) $b:B\rightarrow F(X)$ such that
$x(b) \neq y(b)$.  But then $x \in Z(x(b))$ and $y \in Z(y(b))$ give open neighborhoods that are disjoint. $\Box$
\medskip


\begin{thm}  \label{thm: compact basis}
If $F:{\cal E}\rightarrow {\cal B}$ is a row-finite, strongly surjective KPf and, moreover, every slice category ${\cal B}/ B$ is countable, then for each $\beta\in\E$, 
\begin{enumerate}
\item $Z(\beta)$ is compact and
\item  $F^\infty$ is totally disconnected.
\end{enumerate}
\end{thm}

\noindent{\sc Proof.} Item~(1) implies Item~(2), because $F^\infty$ is Hausdorff and so  $Z(\beta)$ compact implies that it is closed.  Since the ${\cal B}/ B$ are each countable, each $Z(\beta)$ is second countable, so it suffices to show that each $Z(\beta)$ is sequentially compact.

Fix a morphism $\beta:Y\rightarrow X$ in $\cal E$, and order the objects of ${\cal B}/ F(X)$
so as to give a sequence $\{b_i:B_i\rightarrow F(X)\}_{i=0}^\infty$ with $b_0 = F(\beta)$.

Given a sequence $\{x_n\}_{n=0}^\infty$ in $Z(\beta)$, we will construct a convergent subsequence by a diagonalization argument.  Let $\{x_{n,0}\}_{n=0}^\infty := \{x_n\}_{n=0}^\infty$ and observe that 

\[ \forall n\; x_{n,0}(b_0) = \beta, \]

\noindent that is $x_{n,0}(b_0)$ is constant, and thus {\em a fortiori} eventually constant.

Now, suppose we have constructed sequences $\{x_{n,i}\}_{n=0}^\infty$ for $i = 0,\ldots, k$
such that

\begin{enumerate}
\item $\{x_{n,j+1}\}_{n=0}^\infty$ is a subsequence of $\{x_{n,j}\}_{n=0}^\infty$ and
\item $x_{n,i}(b_j)$ is constant for all $j \leq i$.
\end{enumerate}

\noindent Note that the nesting of subsequences ensures that $j \leq i < h\leq k$ implies
$x_{n,i}(b_j) =x_{n,h}(b_j)$.

We construct a subsequence $x_{n,k+1}$ of $x_{n,k}$ satisfying Item (2).  Consider the list
of objects in ${\cal E}\rightarrow X$ given by $\{x_{n,k}(b_{k+1})\}_{n=0}^\infty$.  Each element of this list is a preimage of $b_{k+1}$, but by row-finiteness there are only finitely many distinct preimages. Thus by the pigeonhole principle one, say $x_{N,k}(b_{k+1})$ occurs infinitely often. Let $\{x_{n,k+1}\}_{n=0}^\infty$ be the subsequence of 
all $x_{n,k}$'s for which $x_{n,k}(b_{k+1}) = x_{N,k}(b_{k+1})$.  Now being a subsequence of $\{x_{n,k}\}_{n=0}^\infty$, our new sequence is constant on the $b_j$ for $j \leq k$, and by construction it is constant on $b_{k+1}$, thus $x_{n, k+1}$ satisfies Item (2) as desired.

We claim the diagonal sequence $\{x_{n,n}\}$, necessarily a subsequence of our original 
sequence $\{x_n\}_{n=0}^\infty$, is convergent.  By construction, for each $b:B\rightarrow F(X)$, $x_{n,n}(b)$ is eventually constant, and by Lemma \ref{objectssuffice} for any morphism $(c,d)$ in ${\cal B}/ F(X)$, $x_{n,n}(c,d)$ is eventually constant.

We define the limiting infinite path $x$ by $x(\omega) = \lim_{n\rightarrow \infty} x_{n,n}(\omega)$, whether $\omega$ is an object or morphism of ${\cal B}/ F(X)$, and the limit exists by eventual constancy.  To see $x$ is a functor  let $(a,b),(ab,c)\in \B/ F(X)$.  Then there exists $n$ sufficiently large so that 
\[
x(a,b)x(ab,c)=x_n(a,b)x_n(ab,c)=x_n(a,bc)=x(a,bc);
\]
that is $x$ is a functor.  Likewise that $x$ is a
 section of $F$ follows from its agreement with $x_{n,n}$ for $n$ sufficiently large on each object and morphism. 

It remains only to show that $x = \lim_{n\rightarrow \infty} x_{n,n}$ in the topology on $F^\infty$.  But the basic opens containing $x$ are given by $\{ Z(x(b_i)) : i=0, 1, 2, \ldots \}$, and by construction $x_{n,n} \in Z(x(b_i))$ for all
$n \geq i$, so we are done.  $\Box$ \medskip

Now, having equipped $F^\infty$ with a topology, we consider the continuity of our two families of maps $\res_\mu$ and $\ind_\mu$.

\begin{proposition}
Let $F:\E\to \B$ be a KPf and $\mu\in \E$.  Then $\res_\mu: Z(\mu)\to  Z(s(\mu))$ and $\ind_\mu: Z(s(\mu))\to Z(\mu)$ are continuous.
\end{proposition}

\noindent{\sc Proof.} To see $\res_\mu$ is continuous, it suffices to show that $\res_\mu^{-1}( Z(\nu) \cap Z(s(\mu)))$, the inverse image of generic basic open set for the subspace topology, is open.  

But $Z(\nu) \cap Z(s(\mu)) = \emptyset$ unless $r(\nu)=s(\mu)$, in which
case $Z(\nu) \cap Z(s(\mu)) = Z(\nu)$.  Now $\res_\mu^{-1}( Z(\nu))=\ind_\mu(Z(\nu))=Z(\mu\nu)$ by Proposition~\ref{restrictandinduce} and  is thus open.

To see $\ind_\mu$ is continuous it again suffices to show that $\ind_\mu\inv(Z(\mu)\cap Z(\gamma))$ is open for any $\gamma\in \E$.  Now $Z(\mu)\cap Z(\gamma)$ is either empty or a disjoint union of sets of the form $Z(\mu\nu)$ by Proposition~\ref{disjoint union}, so it suffices to show $\ind_\mu\inv(Z(\mu\nu))$ is open.  But by Proposition~\ref{restrictandinduce}, $\ind_\mu\inv(Z(\mu\nu))=\res_\mu(Z(\mu\nu))=Z(\nu)$ is open as desired.   $\Box$
\medskip

With this in hand we begin our construction of the groupoid.  Let

\[
\mathcal{G}_F:=\{(\mu,\nu, x)\in \E\times \E\times F^\infty: x\in Z(\nu)\quad\text{and}\quad s(\mu)=s(\nu)\}.
\]
We think of $(\mu,\nu,x)$ as a map taking $x$ to $\ind_{\mu}(\res_\nu)(x)$.
We define a relation on $\mathcal{G}_F$ by $(\mu,\nu,x)\sim (\mu',\nu', x')$ if 
\begin{enumerate}
\item $x=x'$, \label{rel: path eq}
\item there exists $\lambda\in \E$ with $x\in Z(\lambda)\subset Z(\nu)\cap Z(\nu')$ and $\ind_\mu \circ \res_\nu|_{Z(\lambda)}=\ind_{\mu'} \circ \res_{\nu'}|_{Z(\lambda)}$,\label{rel: germ eq}
\item there exist $a,b\in \B$ with $F(\mu)a=F(\mu')b$ and $F(\nu)a=F(\nu')b$. \label{rel: base}
\end{enumerate}

A quick check shows that $\sim$ is reflexive, symmetric and transitive so it is an equivalence relation.  Define $G_F=\mathcal{G}_F/\sim$ and denote the image of $(\mu,\nu,x)$ in $G_F$ by $[\mu,\nu,x]$. 

We define composition of morphisms in $G_F$.  As we are thinking of $[\mu,\nu,x]$ as a morphism taking  $x$ to $\ind_{\mu}(\res_\nu)(x)$ we would like to compose  $[\mu,\nu,x],[\sigma,\tau,y]\in G_F\times G_F$ if $x=\ind_\sigma \circ \res_\tau(y)$ and the map should be $y\to \ind_{\mu}\circ \res_{\nu}\circ \ind_\sigma\circ\res_\tau (y)$.  But to define composition in this way we need to find $(\xi,\zeta)\in \E\times \E$ with the map $\ind_{\xi}\circ\res_\zeta=\ind_{\mu}\circ \res_{\nu}\circ \ind_\sigma\circ\res_\tau$ on a neighborhood of $y$ and then define $[\mu,\nu,x]\cdot[\sigma,\tau,y]=[\xi,\zeta,y]$. 

Consider $[\mu,\nu,x],[\sigma,\tau,y]\in G_F\times G_F$ with $x=\ind_\sigma \circ \res_\tau(y)$. We construct a candidate for $[\xi,\zeta,y]$.   Since $x=\ind_\sigma \circ \res_\tau(y)$,  $x(F(\nu))=\nu$ and $x(F(\sigma))=\sigma$.  Thus $F(\nu), F(\sigma)$ is a cospan in the right Ore category $\B$ so there exists $a,b\in \B$ with $F(\nu)a=F(\sigma)b$.  Thus $\nu x_2(F(\nu,a))=x(F(\nu)a)=x(F(\sigma)b)=\sigma x_2(F(\sigma),b)$.  Take $\gamma=x_2(F(\nu,a))$ and $\eta=x_2(F(\sigma),b)$ so that $\nu\gamma=\sigma\eta$. 

Notice that on $Z(\tau\eta)$ we have

\begin{align*}
\ind_{\mu}\circ \res_\nu\circ\ind_{\sigma}\circ\res_{\tau}&=\ind_{\mu}\circ \ind_{\gamma}\circ \res_\gamma \circ\res_\nu\circ\ind_{\sigma}\circ \ind_{\eta}\circ \res_\eta \circ\res_{\tau}\\
&=\ind_{\mu\gamma}\circ \res_{ \nu\gamma}\circ\ind_{\sigma\eta}\circ \res_{\tau\eta}\\
&=\ind_{\mu\gamma}\circ \res_{\tau\eta}.
\end{align*}

\begin{lem}
\label{lem: mult def}
The formula 

\[
[\mu,\nu,\ind_\sigma \circ \res_\tau(y)][\sigma,\tau,y] := [\mu\gamma,\tau\eta,y],
\]

\noindent for $\eta$ and $\gamma$ as in the discussion above, with $x = \ind_\sigma \circ \res_\tau(y)$, gives a well-defined composition in $G_F$.
\end{lem}

\noindent {\sc Proof.} Suppose $[\sigma,\tau,y]=[\sigma',\tau',y]$ and $[\mu,\nu,\ind_\sigma \circ \res_\tau(y)]=[\mu',\nu',\ind_{\sigma'} \circ \res_{\tau'}(y)]$.  Chose $\gamma',\eta'$ such that $\nu'\gamma'=\sigma'\eta'$.  Then 
\[
[\mu,\nu,\ind_\sigma \circ \res_\tau(y)][\sigma,\tau,y]=[\mu\gamma,\tau\eta,y]\quad\text{ and}\quad [\mu',\nu',\ind_{\sigma'} \circ \res_{\tau'}(y)][\sigma',\tau',y]=[\mu'\gamma',\tau'\eta',y],
\]
we must show $[\mu\gamma,\tau\eta,y]=[\mu'\gamma',\tau'\eta',y]$.  We have $y=y$ and, since composition of germs is well-defined, there exists a $\lambda\in \E$ so that on $Z(\lambda)$, $\ind_{\mu\gamma}\circ \res_{\tau\eta}=\ind_{\mu}\circ \res_\nu\circ\ind_{\sigma}\circ\res_{\tau}=\ind_{\mu'}\circ \res_{\nu'}\circ\ind_{\sigma'}\circ\res_{\tau'}=\ind_{\mu'\gamma'}\circ \res_{\tau'\eta'}$.  It remains to show there exist $R,R'$ such that $F(\mu\gamma)R=F(\mu'\gamma')R'$ and $F(\tau\eta)R=F(\tau'\eta')R'$.  Let

\begin{align*}
m&=F(\mu)& m'&=F(\mu')\\
n&=F(\nu)& n'&=F(\nu')\\
w&=F(\sigma) &w'&=F(\sigma')\\
t&=F(\tau) & t' &=F(\tau')\\
g&=F(\gamma) & g'&=F(\gamma')\\
h&=F(\eta) & h'&=F(\eta').\\
\end{align*}
So that
\begin{align}
\label{1} ng&=wh\\
\label{2}n'g'&=w'h'.
\end{align}
 Since $[\mu,\nu, \ind_\sigma\circ \res_\tau(y)]=[\mu',\nu',\ind_{\sigma'}\circ \res_{\tau'}(y)]$ and $[\sigma,\tau,y ]=[\sigma',\tau',y]$ there exist $p,p',q,q'\in\B$ such that 
\begin{align}
\label{3} mp&=m'p' & np&=n'p'\\
\label{4} wq&=w'q' & tq &=t'q'.
\end{align} 

\noindent Since $(p,g), (p',g'), (q,b),(q',b')$ are cospans, by the right Ore condition, there exist 
\[a_1, a_2, b_1,b_2, a_1' a_2', b_1', b_2'\in \B \text{  such that}\]

\begin{align}
\label{5} pa_1 &= gb_1 & qa_2 &= hb_2\\
\label{6} p'a_1' &= g'b_1' & q'a_2' &= h'b_2'.
\end{align}

Since $s(p)=s(p'), s(q)=s(q')$, we have $(a_1,a_1'), (a_2,a_2')$ are cospans. Using the right Ore condition again, there exist $c_1, c_1', d_1, d_1'\in \B$ such that
\begin{align}
\label{7} a_1 c_1&=a_1' c_1' & a_2 d_1&=a_2'd_1'.
\end{align}

Since $r(b_1)=s(g)=s(h)=r(b_2)$ and $r(b_1')=s(g')=s(h')=r(b_2')$, we have $(b_1c_1, b_2d_1)$ and $(b_1'c_1', b_2d_1')$ are cospans. Using the right Ore condition, there exist $c_2, c_2', d_2, d_2'\in \B$ such that
\begin{align}
\label{8} b_1 c_1 c_2&=b_2 d_1 d_2 & b_1' c_1' c_2'&=b_2' d_1' d_2'.
\end{align}

Define

\begin{align*}
M&=ngb_1c_1=npa_1c_1=n'p'a_1'c_1'=n'g'b_1'c_1', \quad\text{and}\\
N&=whb_2d_1=wqa_2d_1=w'q'a_2'd_1'=w'b'b_2'd_1'.\\
\end{align*}
Now $M,N$ is a cospan and $Mc_2=Nd_2$ and $Mc_2'=Nd_2'$.  Thus by strongly right Ore there exist $k,\ell\in \B$ with 

\begin{align}
\label{9}c_2k&=c_2'\ell & d_2k &=d_2'\ell.
\end{align}
Take 
\begin{align*}
R&=b_1c_1c_2k & R'&=b_1'c_1'c_2'\ell.
\end{align*}
Then
\begin{align*}
mgR &=mgb_1c_1c_2k=mpa_1c_1c_2k \quad\text{by \eqref{5}}\\
&=mpa_1c_1c_2'\ell\quad\text{by \eqref{9}}\\
&=mpa_1'c_1'c_2'\ell \quad\text{by \eqref{7}}\\
&=m'p'a_1'c_1'c_2'\ell\quad\text{by \eqref{3}}\\
&=m'g'b_1'c_1'c_2'\ell\quad\text{by \eqref{6}}\\
&=m'g'R'
\end{align*}
and
\begin{align*}
thR &=thb_1c_1c_2k =thb_2d_1d_2k \quad\text{by \eqref{8}}\\
&=tqa_2d_1d_2k\quad\text{by \eqref{5}}\\
&=tqa_2d_1d_2'\ell\quad\text{by \eqref{9}}\\
&=tqa_2'd_1'd_2'\ell\quad\text{by \eqref{7}}\\
&=t'q'a_2'd_1'd_2'\ell\quad\text{by \eqref{4}}\\
&=t'h'b_2'd_1'd_2'\ell\quad\text{by \eqref{6}}\\
&=t'h'b_1'c_1'c_2'\ell\quad\text{by \eqref{8}}\\
&=t'h'R'.
\end{align*}
This gives that multiplication is well-defined on $G_F$.$\Box$ \medskip

\begin{lem}
Under the multiplication defined in Lemma~\ref{lem: mult def}, $G_F$ is a groupoid.
\end{lem}
\noindent {\sc Proof.} 
Notice for $[\mu,\nu,x]\in G_F$, $[r(\mu), r(\mu),\ind_\mu(\res_\nu(x))][\mu,\nu,x]=[r(\mu)\mu,\nu s(\mu),x]=[\mu,\nu,x]$ and $[\mu,\nu,x][r(\nu),r(\nu),x]=[\mu s(\nu),r(\nu)\nu,x]=[\mu,\nu,x]$;  that is elements of the form $[r(x),r(x),x]$ act as units in $G_F$ and the map $x\mapsto [r(x),r(x),x]$ identifies $F^\infty$ with the unit space of $G_F$.    Also $[r(x),r(x),x]=[x(a),x(a),x]$ for any object $a$ in $ \B/ F(r(x))$.

Now for $[\mu,\nu,x]\in G_F$ consider $[\nu,\mu, \ind_\mu(\res_\nu(x))]$, then
\begin{align*}
[\mu,\nu,x][\nu,\mu, \ind_\mu(\res_\nu(x))]&=[\mu,\mu, \ind_\mu(\res_\nu(x))]\quad\text{and}\\
[\nu,\mu, \ind_\mu(\res_\nu(x))][\mu,\nu,x]&=[\nu,\nu,x]
\end{align*}
that is the inverse of $[\mu,\nu,x]$ is $[\nu,\mu \ind_\mu(\res_\nu(x))]$.
Thus $G_F$ is a groupoid. $\Box$ \medskip

\begin{rmk}
Let $F:\E\to \B$ be a row-finite Kumjian-Pask fibration where $\B$ (and hence $\E$) is left and right cancellative. Suppose that $\E$ has no inverses, then $\E$ is a  finitely aligned category of paths in the sense of \cite[Definition~3.1]{S}. Let $X$ be an object in $\E$.  Define $A_X$ to be the set of finite disjoint unions of sets of the form $(\bigcup_{j=1}^N \alpha_j\E)-(\bigcup_{k=1}^M\beta_k\E)$ and $\Omega_X$ be the set of ultra filters on $A_X$: that is $\omega\in \Omega_X$ is a subset of the power set on $A_X$ with the property that for every $E\in A_X$ either $E\in \omega$ or there is an $F\in  \omega$ such that  $F\cap E=\emptyset$.  Take $\Omega=\cup_{X\in \text{Obj}(\E)} \Omega_X$.  For $\beta \in X\E$ we can define a set $\hat{\beta\E}=\{\omega\in\Omega_X: \beta\E\in \omega\}$.  The set  $\{\hat{\beta\E}\}_{\beta\in \E}$ forms a subbasis for a topology on $\Omega$. For $\alpha\in \E$, there is a map $\tilde{\alpha}: A_{s(\alpha)}\to A_{r(\alpha)}$ characterized by $\tilde{\alpha}(\beta\E)=(\alpha\beta)\E$; $\tilde{\alpha}$ then induces a continuous map $\hat{\alpha}: \Omega_{s(\alpha)}\to \Omega_{r(\alpha)}$.   In \cite{S}, Spielberg defines a groupoid $G(\E)$ to be the groupoid generated by the germs of the maps $\tilde{\alpha}$ for all $\alpha\in \E$.  The unit space of $G(\E)$ can then be identified with $\Omega$.

We can view $G_F^{(0)}=F^\infty\subset \Omega=G(\E)^{(0)}$ by taking $x\in F^\infty$ to 
\[
\omega_x=\{E\in A_{r(x)}: x(a)\E\subset E\quad\text{for some}\quad a\in \B, r(a)=F(r(x))\}.
\]
To see that $\omega_x$ is an ultrafilter on $A_{r(x)}$ we need to see that if $E\in A_{r(x)}$ then either there exists $a$ with $x(a)\E\subset E$ or $x(a)\E\cap E=\emptyset$.  Suppose $x(b)\E\cap E\neq \emptyset$ for all $b\in \B$ with $r(b)=F(r(x))$.  Now $E=\bigcup_{i=1}^Q( (\bigcup_{j=1}^{N_i} \alpha_{i,j}\E)-(\bigcup_{k=1}^{M_i}\beta_{i,k}\E))$.  Since $\B$ is right Ore, there exists $a\in \B$ that extends $F(\alpha_{i,j})$ and $F(\beta_{i,k})$ for all choices of $i,j,$ and $k$.  We claim $x(a)\E\subset E$.  We know that $x(a)\E\cap E\neq \emptyset$.  Thus there exists $i_0,j_0$ sch that $x(a)\E\cap (\alpha_{i_0,j_0}\E-(\bigcup_{k=1}^{M_{i_0}}\beta_{i_0,k}\E))\neq \emptyset$.  Let $\gamma\in x(a)\E\cap (\alpha_{i_0,j_0}\E-(\bigcup_{k=1}^{M_{i_0}}\beta_{i_0,k}\E))$ so that $\gamma=x(a)\gamma'=\alpha_{i_0,j_0}\gamma''$.  Since $a$ extends $F(\alpha_{i_0,j_0})$, by unique factorization we get $x(a)=\alpha_{i_0,j_0}\eta$ for some $\eta$.  Thus $x(a)\E\subset \alpha_{i_0,j_0}\E$.  Now $\gamma\not\in \beta_{i_0,k}\E$ for all $k$.  Since $a$ extends $F(\beta_{i_0,k})$ for all $k$ we must have $x(F(\beta_{i_0,k}))\neq \beta_{i_0,k}$ for all $k$ so that $x(a)\E \cap \beta_{i_0,k}\E=\emptyset$ for all $k$.  That is $x(a)\E\subset\alpha_{i_0,j_0}\E-(\bigcup_{k=1}^{M_{i_0}}\beta_{i_0,k}\E)\subset E$ as desired.  So that $\omega_x$ is an ultra filter and we can view $F^\infty\subset \Omega$ as claimed.
\end{rmk}
 
We want to define a topology on $G_F$.  For $\mu,\nu\in \E$ consider the set
\[
Z(\mu,\nu):=\{[\alpha,\beta, x]\in G_F: [\alpha,\beta,x]=[\mu,\nu,x]\}.
\]
That is $[\alpha,\beta,x]\in Z(\mu,\nu)$ if $x\in Z(\nu)$ and $(\alpha,\beta,x)\sim (\mu,\nu,x)$.

Given $(\alpha,\beta), (\mu,\nu)\in \E\times\E$, notice that if there exists an $a$ such that $F(\alpha)a=F(\mu)$ and $F(\beta)a=F(\nu)$, then $Z(\mu,\nu)\subset Z(\alpha,\beta)$ if there exists a $\gamma$ with $F(\gamma)=a$ and $(\alpha\gamma,\beta\gamma)=(\mu,\nu)$, and $ Z(\mu,\nu)\cap Z(\alpha,\beta)=\emptyset$, otherwise.

Indeed if $Z(\mu,\nu)\cap Z(\alpha,\beta)\neq\emptyset$ then there exists $[\mu,\nu,x]\in Z(\alpha,\beta)$. So $x(F(\nu))=x(F(\beta)a)=x(F(\beta))x_2(F(\beta),a)$; so by unique factorization we have $x(F(\beta))=\beta$ and $\nu=\beta \gamma$ with $\gamma=x_2(F(\beta),a)$. Now 
\[
\mu=\ind_{\alpha}\circ \res_{\beta} (x)(F(\alpha)b)=\alpha \res_\beta(x)(b)=\alpha\gamma.
\]

\begin{proposition}
\label{disjoint basis}
Suppose that $Z(\alpha,\beta)\cap Z(\sigma,\tau)\neq \emptyset$.  Then there exist $a,b\in \B$ and $I=\{(\alpha\gamma,\beta\gamma)\in \E\times \E: F(\gamma)=a\}\cap \{(\sigma\eta,\tau\eta)\in \E\times \E: F(\eta)=b\}\neq \emptyset$ such that 

\[
Z(\alpha,\beta)\cap Z(\sigma,\tau)=\bigcup_{(\mu,\nu)\in I} Z(\mu,\nu)
\]
and the union on the right hand side is disjoint.
\end{proposition}
 
\noindent {\sc Proof.}
Since $Z(\alpha,\beta)\cap Z(\sigma,\tau)\neq \emptyset$ there exists $x\in F^\infty$ with $[\alpha,\beta,x]=[\sigma,\tau,x]$.  This occurs if and only if $\ind_\alpha\circ \res_\beta=\ind_\sigma\circ\res_\tau$ on some neighborhood and there exist $a,b\in \B$ such that 
\[ F(\alpha)a=F(\sigma) b\quad\text{and}\quad F(\beta)a=F(\tau)b.\]
Take $I$ for this $a,b\in \B$. By the definition of $I$ if $(\mu,\nu)\in I$ then $Z(\mu,\nu)\subset Z(\alpha,\beta)\cap Z(\sigma,\tau)$ and so $\bigcup_{(\mu,\nu)\in I} Z(\mu,\nu)\subset Z(\alpha,\beta)\cap Z(\sigma,\tau)$.

Now if we assume $[\alpha,\beta,x]=[\sigma,\tau,x]\in Z(\alpha,\beta)\cap Z(\sigma,\tau)\neq \emptyset$.  Take $\gamma=x_2(F(\beta),a)$ and $\eta=x_2(F(\tau),b)$.  By definition we then have $\beta\gamma=x(F(\beta)a)=x(F(\tau)b)=\tau\eta$, and since $\ind_\alpha\circ \res_\beta=\ind_\sigma\circ\res_\tau$ on a neighborhood of $x$ we also have 
\begin{align*}
\alpha\gamma&=\alpha x_2(F(\beta),a)\\
&=\alpha\res_\beta(x)(a)\\
&=\ind_\alpha\circ\res_\beta(x)(F(\alpha)a)\\
&=\ind_\sigma\circ\res_\tau(x)(F(\sigma)b)\\
&=\sigma\eta.
\end{align*}
That is $(\alpha\gamma,\beta\gamma)=(\sigma\eta,\tau\eta)\in I$ and so $[\alpha,\beta,x]\in \bigcup_{(\mu,\nu)\in I} Z(\mu,\nu)$. Thus $Z(\alpha,\beta)\cap Z(\sigma,\tau)=\bigcup_{(\mu,\nu)\in I} Z(\mu,\nu)$ as desired.  Now for $(\mu,\nu),(\mu',\nu')\in I$, $F(\mu)=F(\mu'), F(\nu)=F(\nu')$ so that $Z(\mu,\nu)\cap Z(\mu',\nu')=\emptyset$ if $(\mu,\nu)\neq(\mu',\nu')$ by unique factorization.
$\Box$ \medskip

\begin{cor}
\label{lem: top}
The set $\{Z(\mu,\nu)\}$ is  a basis for a topology on $G_F$. 
\end{cor}

\begin{lem}
\label{lem: comp cont}
With respect to this topology composition is continuous and inversion is a homeomorphism: that is $G_F$ is a topological groupoid. 
\end{lem}

\noindent {\sc Proof.}
To see inversion is continuous notice that $Z(\mu,\nu)\inv=Z(\nu,\mu)$, since inversion is an involution it is a homeomorphism.

To see composition is continuous, let $Z(\mu,\nu)$ be a basic open set and $[\alpha,\beta,\ind_\sigma(\res_\tau(x))][\sigma,\tau,x]=[\alpha\gamma,\tau\eta,x]\in Z(\mu,\nu)$; in particular $\beta\gamma=\sigma\eta$. We need to find neighborhoods $Z(\alpha',\beta')$ and $Z(\sigma',\tau')$ of $[\alpha,\beta,\ind_\sigma(\res_\tau(x))]$ and $[\sigma,\tau,x]$ respectively so that $Z(\alpha',\beta')Z(\sigma',\tau')\subset Z(\mu,\nu)$.

Since $[\alpha\gamma,\tau\eta,x]\in Z(\mu,\nu)$, there exist $d,e\in \B$ so that 
\begin{align*}
F(\alpha\gamma)d&=F(\mu)e\quad\text{and}\\
F(\tau\eta)d &=F(\nu)e.
\end{align*}
Further there exists $k\in \B$ such that 
\[\ind_\mu\circ\res_\nu|_{Z(x(F(\nu)ek))}=\ind_{\alpha\gamma}\circ\res_{\tau\eta}|_{Z(x(F(\tau\eta)dk))}.
\]
Chose $\delta,\epsilon, \kappa\in \E$ so that 
\[
\nu\epsilon\kappa=x(F(\nu)ek)=x(F(\tau\eta)dk)=\tau\eta\delta\kappa
\]
and $F(\delta)=d,F(\epsilon)=e, F(\kappa)=k$.  Since $\beta\gamma=\sigma\eta$, $s(\gamma)=s(\eta)=r(\delta)$.  Take 
\[
\alpha'=\alpha\gamma\delta\kappa,\quad\quad \beta'=\beta\gamma\delta\kappa,\quad\quad\sigma'=\sigma\eta\delta\kappa,\quad\text{and} \quad\quad\tau'=\tau\eta\delta\kappa.
\]
By construction $\beta'=\sigma'$.  So by definition 
\[
Z(\alpha',\beta')Z(\sigma',\tau')=Z(\alpha',\tau')=Z(\alpha\gamma\delta\kappa,\tau\eta\delta\kappa)=Z(\alpha\gamma\delta\kappa,\nu \epsilon\kappa).
\]
We what to show $Z(\alpha',\beta')Z(\sigma',\tau')\subset Z(\mu,\nu)$ and by the above computation it suffices to show $\alpha\gamma\delta\kappa=\mu\epsilon\kappa$.  But 
\begin{align*}\alpha\gamma\delta\kappa&=\ind_{\alpha\gamma}\circ\res_{\tau\eta}(x)(F(\alpha\gamma)dk)=\ind_\mu\circ\ind_\nu(x)(F(\alpha\gamma)dk)\\
&=\ind_\mu\circ\ind_\nu(x)(F(\mu)ek)=\mu\epsilon\kappa\
\end{align*}
as desired. 
$\Box$ \medskip
 
\begin{lem}
\label{lem: etale}
For $\mu,\nu\in \E$.  The maps $r|_{Z(\mu,\nu)}: Z(\mu,\nu)\to Z(\mu)$ and $s|_{Z(\mu, \nu)}: Z(\mu,\nu)\to Z(\nu)$ are homeomorphisms and $G_F$ is \'etale.
\end{lem}   

\noindent {\sc Proof.}
We begin by showing $s|_{Z(\mu, \nu)}: Z(\mu,\nu)\to Z(\nu)$ is a homeomorphism.  First note that $s|_{Z(\mu,\nu)}$ is injective.  Now since $Z(\alpha,\beta)$ form a basis for $G_F$ (and therefore for $Z(\mu,\nu)$) and $Z(\beta)$ form a basis for $F^\infty$ (and hence for $Z(\nu)$) it suffices to see $s(Z(\alpha,\beta))=Z(\beta)$.  Now if $x\in Z(\beta)$  then $[\alpha,\beta,x]\in Z(\alpha,\beta)$ and so $Z(\beta)\subset s(Z(\alpha,\beta))$.  By definition $s(Z(\alpha,\beta))\subset Z(\beta)$, so that $s(Z(\alpha,\beta))=Z(\beta)$ as desired.  Thus $s|_{Z(\mu, \nu)}: Z(\mu,\nu)\to Z(\nu)$ is a homeomorphism.  Since $r([\alpha,\beta,x])=s([\alpha,\beta,x]\inv)$ and inversion is a homeomorphism we get that $r|_{Z(\mu,\nu)}: Z(\mu,\nu)\to Z(\mu)$ is a homeomorphism as well. $G_F$ is now \'etale by definition.
$\Box$ \medskip

\begin{thm}Let $F:\E\to \B$ be a row-finite Kumjian-Pask fibration with every slice category $\B/B$ countable.  Then $G_F$ with the topology induced by the basis $\{ Z(\mu, \nu) : \mbox{$(\mu, \nu)$ is a span in } {\cal E}\} $ is totally disconnected locally compact Hausdorff.
\end{thm}

\noindent{\sc Proof.}  Note that $Z(\mu,\nu)$ is homeomrphic to $Z(\nu)$ by Lemma~\ref{lem: etale} and $Z(\nu)$ is compact from Theorem~\ref{thm: compact basis}.  Thus $G_F$ is locally compact; it will also follow that $G_F$ is totally disconnected once we show that $G_F$ is Hausdorff.  To show that $G_F$ is Hausdorff, first note that if two morphisms $[\alpha, \beta, x]$ and $[\gamma, \delta, x]$ both lie in a basic open $Z(\mu, \nu)$ then they are equal, since both are $[\mu, \nu, x]$.   

For brevity, having chosen to denote morphisms in $\cal E$ with Greek letters, we denote their image under $F$ with the corresponding Latin letter (thus, for instance $F(\alpha)$ will be denoted $a$).

Now, suppose $[\alpha, \beta, x] \neq [\gamma, \delta, y]$.  First suppose $x\neq y$. If $r(\beta) \neq r(\delta)$ or $r(\alpha) \neq r(\delta)$, it is immediate that $Z(\alpha, \beta)$ and $Z(\gamma, \delta)$ are disjoint.  So suppose $\alpha, \beta$ and $\gamma, \delta$ are two spans between a pair of objects $\Gamma = r(\gamma) = r(\alpha)$ and $\Delta = r(\delta) = r(\beta)$ in $\cal E$.  Since $x \neq y$, there exists a morphism $f:X\rightarrow D$ in $\cal B$ for which $x(f) \neq y(f)$ are distinct lifts. 

By the right Ore condition, we can complete the cospan $b, d$ to a commutative square $bk = d\ell$,
and the cospan $bk = d\ell, f$ to a commutative square $fn = bkm (= d\ell m)$.  By the unique factorization lifting property, it follows that $x(f) \neq y(f)$ implies $x(bkm) = x(fn) \neq y(fn) = y(d\ell m)$.
But $[\alpha, \beta, x] = [\alpha x(km), x(fn), x] \in Z(\alpha x(km), x(fn))$ and $[\gamma, \delta, y] = 
[\gamma y(\ell m), y(fn), y] \in Z(\gamma y(\ell m), y(fn))$, and these two open sets are disjoint, since all
elements of the first have infinite paths lifting $fn$ to $x(fn)$, while all elements of the second have infinite paths lifting $fn$ to $y(fn)$.

Now assume that $x=y$. If $Z(\alpha, \beta)$ and $Z(\gamma, \delta)$ are disjoint we have separated the two morphisms with disjoint opens.  If not, there exists a 
$[\epsilon, \zeta, z] \in Z(\alpha, \beta) \cap Z(\gamma, \delta)$, in which case 

\[ [\alpha, \beta, z] = [\gamma, \delta, z] = [\epsilon, \zeta, z]. \]
\noindent From the equality of $ [\alpha, \beta, z]$ and $[\gamma, \delta, z]$, we have a span in $\cal B$
from $s(a) = s(b)$ to $s(c) = s(d)$, $m, n$ such that $am = cn$ and $bm = dn$.  

The square  $bm = dn$ admits a lift by $x$ (as a span in the slice category ${\cal B} /
r(m)$).  Let $\omega = x(bm) = x(dm)$.  Now $(\ind_\alpha\circ\res_\beta(x))(am)=\alpha x_2(b,m)$ and $(\ind_\gamma\circ\delta(x))(cn)=\gamma x_2(d,n)$. Now  $\alpha x_2(b,m)\neq \gamma x_2(d,n)$ otherwise $\ind_\alpha \res_\beta = \ind_\gamma \res_\delta$ on $Z(\omega)$ and we get $[\alpha, \beta, x] = [\gamma, \delta, x]$.

Note, then that $\alpha x_2(b,m)$ and $ \gamma x_2(d,n)$ are then two unequal lifts of $am = cn$.  It is clear that
$[\alpha, \beta, x] = [\alpha x_2(b,m), \omega, x] \in Z(\alpha x_2(b,m), \omega)$ and $[\gamma, \delta, x] = [\gamma x_2(d,n), \omega, x] \in Z(\gamma x_2(d,n), \omega)$.  But $Z(\alpha x_2(b,m), \omega)$ and $Z(\gamma x_2(d,n), \omega)$ are disjoint, since the image of any element of $Z(\omega)$ under any element of the first lifts $am = cn$ to
$\alpha x_2(b,m)$, while the images of elements of $Z(\omega)$ under any element of the latter lifts $am = cn$ to $\gamma x_2(d,n)$ by Proposition~\ref{restrictandinduce}.
$\Box$
\medskip

\begin{example}
If $d:\Lambda\to \N^k$ is a $k$-graph, Kumjian and Pask \cite{KP} define $G_\Lambda:=\{(x,m-n,y)\in d^\infty\times \Z\times d^\infty: \res_{x(m)}(x)=\res_{y(n)}(y)\}$ and the topology on $G_\Lambda$ is generated by $\mathcal{Z}(\alpha,\beta):=\{(\ind_\alpha z,d(\alpha)-d(\beta),\ind_\beta y): r(y)=s(\beta)\}$.  This on the surface is different from $G_{d}$ as defined above.  However the map $[\alpha,\beta,x]\mapsto (\ind_\alpha\circ \res_\beta(x), d(\alpha)-d(\beta),x)$ defines an isomorphism of $G_d$  to Kumjian and Pask's $G_\Lambda$.  This map is well-defined since $[\alpha,\beta,x]=[\mu,\nu,x]$ gives that $\ind_\alpha\circ\res_\beta=\ind_\mu\circ \res_\nu$ on a neighborhood of $x$ and there exists $m,n$ such that $d(\alpha)+m=d(\mu)+n$ and $d(\beta)+m=d(\nu)+n$, in particular $\ind_\alpha\circ \res_\beta(x)=\ind_\mu\circ \res_\nu(x)$ and $d(\alpha)-d(\beta)=d(\mu)-d(\nu)$; that is $(\ind_\alpha\circ \res_\beta(x), d(\alpha)-d(\beta),x)=(\ind_\mu\circ \res_\nu(x), d(\mu)-d(\nu),x)$.  A quick computation shows that this map is a homomorphism and takes the sets $Z(\alpha,\beta)$ to $\mathcal{Z}(\alpha,\beta)$ so that it is also a homeomorphism.
\end{example}

\begin{example}
Let $G$ be a groupoid then $G$ in particular is a small category with objects $G^{(0)}$.  Then $\Id_G:G\to G$ is a KPf.  For each $x\in G^{(0)}$ there is exactly one section $x: G/ \Id_{G}(x)\to G/ x$.  So we can identify $G^{(0)}$ with $\Id_{G}^\infty$ as sets.  For $\gamma\in G$ the map $\gamma\mapsto [\gamma, s(\gamma),s(\gamma)]$ is a bijective homomorphism of $G$ to $G_{\Id_{G}}$.  Indeed, since $[\alpha,\beta,r(\beta)]=[\alpha\beta\inv,s(\alpha\beta\inv),s(\alpha\beta\inv)]$ we have that the map is onto.  To see it is injective, note that for $[\gamma,s(\gamma),s(\gamma)]=[\eta,s(\eta),s(\eta)]$ there exists $\sigma,\tau\in G$ such that $\Id_{G}(\gamma)\sigma=\Id_G(\eta)\tau$ and $s(\gamma)\sigma=s(\eta)\tau$; that is $\sigma=\tau$ and so $\gamma=\eta$.  It is an homomorphism since $[\gamma,s(\gamma),s(\gamma)][\eta,s(\eta),s(\eta)]=[\gamma\eta,\eta,s(\gamma)][\eta,s(\eta),s(\eta)]=[\gamma\eta,s(\eta),s(\eta)]$.  

Thus $G$ is isomorphic to $G_{\Id_G}$ as groupoids. But this isomorphism is not necessarily of topological groupoids.  Indeed we know that $G_{\Id_G}$ is totally disconnected, so it can be isomorphic to $G$ only if $G$ is totally disconnected.  But digging a bit deeper, since each object in $G$ corresponds to exactly one infinite path, we have $Z(\alpha)=\{r(\alpha)\}$ for all elements $\alpha\in G$.  Therefore the unit space in $G_{\Id_G}$ has the discrete topology and therefore $G_{\Id_G}$ has the discrete topology.  Thus we have that $G\cong G_{\Id_G}$ as topological groupoids if and only if $G$ is discrete.  In particular if $H$ is a discrete group $H\cong G_{\Id_H}$.
\end{example}

\section{Isomorphism of $C^*$-algebras}

In this section we show that the $C^*(F)\cong C^*(G_F)$ if $\B$ is left cancellative.

\begin{lem}
\label{hom}
Let $F:\E\to \B$ be a row-finite Kumjain Pask fibration.  Suppose all the morphisms in $\B$ are monic.  Then the map $s_\alpha\mapsto \ch_{Z(\alpha,s(\alpha))}$ extends to a well-defined homomorphism of $\Upsilon: C^*(F)\to C^*(G_F)$
\end{lem}

\noindent {\sc Proof.}
To prove the lemma it suffices to show that $\{\ch_{Z(\alpha,s(\alpha))}\}$ is a Cuntz-Krieger $F$-family in $C^*(G_F)$.  For this we compute.  First note that
\begin{align*}
\ch_{Z(\alpha, s(\alpha))}^*[\mu,\nu,x]&=\ch_{Z(\alpha, s(\alpha))}[\nu,\mu, \ind_{\mu}\circ\res_{\nu}(x)]\\
&={\begin{cases} 1 &\text{if~} [\nu,\mu, \ind_{\mu}\circ\res_{\nu}(x)]=[\alpha,
s(\alpha) \ind_{\mu}\circ\res_{\nu}(x)]\\
0 &\text{otherwise}
\end{cases}}\\
&={\begin{cases} 1 &\text{if~} [\mu,\nu,x]=[s(\alpha),\alpha,x]\\
0 &\text{otherwise}
\end{cases}}\\
&=\ch_{Z(s(\alpha),\alpha)}[\mu,\nu,x].
\end{align*}
Given $a\in \B$, we want to show $\sum_{F(\alpha)=a}\ch_{Z(\alpha, s(\alpha))}\ch_{Z(\alpha,s(\alpha))}^*=\ch_{Z(r(\alpha),r(\alpha))}$. Since $Z(r(\alpha))=Z(r(\alpha), r(\alpha))$ is the disjoint union $\bigcup_{F(\alpha)=a} Z(\alpha,\alpha)$, it suffices to see  that $\ch_{Z(\alpha, s(\alpha))}\ch_{Z(s(\alpha),\alpha)}=\ch_{Z(\alpha,\alpha)}$.  We compute 
\begin{align*}
\ch_{Z(\alpha,s(\alpha))}*&\ch_{Z(s(\alpha),\alpha)}[\xi,\zeta,x]=\sum_{r([\mu,\nu,y])=r([\xi,\zeta,x])}\ch_{Z(\alpha,s(\alpha))}([\mu,\nu,y])\ch_{Z(s(\alpha),\alpha)}([\mu,\nu,y]\inv[\xi,\zeta,x])\\
&={\begin{cases} \ch_{Z(s(\alpha),\alpha)}([\alpha,s(\alpha),\res_\alpha\circ\ind_\xi\circ\res_\zeta x]\inv[\xi,\zeta,x]) &\text{if~~} \ind_\xi\circ\res_\zeta x\in Z(\alpha) \\
0 & \text{otherwise}
\end{cases}}\\
&={\begin{cases} \ch_{Z(s(\alpha),\alpha)}([s(\alpha),\alpha, \ind_\xi\circ\res_\zeta x][\xi,\zeta,x]) &\text{if~~} \ind_\xi\circ\res_\zeta x\in Z(\alpha) \\
0 & \text{otherwise}.
\end{cases}}
\intertext{For $[s(\alpha),\alpha \circ\ind_\xi\circ\res_\zeta x],[\xi,\zeta,x]$ to be composable there exists $\gamma,\eta\in E$ such that $\alpha\gamma=\xi\eta$ and}
&={\begin{cases} \ch_{Z(s(\alpha),\alpha)}([\gamma,\zeta\eta, x]) &\text{if~} \ind_\xi\circ\res_\zeta x\in Z(\alpha) \\
0 & \text{otherwise}.
\end{cases}}
\end{align*}
This is zero unless $[\gamma, \zeta\eta,x]=[s(\alpha),\alpha,x]$; in particular $x\in Z(\alpha)$.  We want to show this implies  $[\xi,\zeta,x]=[\alpha,\alpha,x]$ so that $\ch_{Z(\alpha,s(\alpha))}^**\ch_{Z(\alpha,s(\alpha))}=\ch_{Z(\alpha,\alpha)}$.  Now $[\gamma, \zeta\eta,x]=[s(\alpha),\alpha,x]$ implies there is a suitable neighborhood such that $\ind_\gamma\circ \res_{\zeta\eta}=\res_\alpha$ which happens if and only if $\res_{\zeta\eta}=\res_{\alpha\gamma}$ on this neighborhood.  But we assumed that $\alpha\gamma=\xi\eta$ so $\res_{\zeta\eta}=\res_{\xi\eta}$ which implies $\ind_\xi\circ \res_\zeta=\Id_{Z(\xi)}=\Id_{Z(\alpha)}=\ind_\alpha\circ\res_\alpha$ on this neighborhood.
Now we also know there exists $a,b\in \B$ such that $F(\gamma)a=b$ and $F(\zeta\eta)a=F(\alpha)b=F(\alpha\gamma)a=F(\xi\eta)a$.  Thus we have $[\xi,\zeta,x]=[\alpha,\alpha,x]$ as desired.

It remains to see that $ \ch_{Z(\alpha,s(\alpha))}^**\ch_{Z(\alpha,s(\alpha))}=\ch_{Z(s(\alpha),\alpha)}*\ch_{Z(\alpha,s(\alpha))}=\ch_{Z(s(\alpha),s(\alpha))}$. For this we compute.
\begin{align*}
\ch_{Z(s(\alpha),\alpha)}*&\ch_{Z(\alpha,s(\alpha))}[\xi,\zeta,x]=\sum_{r([\mu,\nu,y])=r([\xi,\zeta,x])}\ch_{Z(s(\alpha),\alpha)}([\mu,\nu,y])\ch_{Z(\alpha,s(\alpha))}([\mu,\nu,y]\inv[\xi,\zeta,x])\\
&={\begin{cases} \ch_{Z(\alpha,s(\alpha))}([s(\alpha),\alpha,\ind_{\xi}\circ \res_{\zeta} x]\inv[\xi,\zeta,x]) &\text{if }\ind_{\xi}\circ \res_{\zeta} x\in Z(s(\alpha))\\
0 &\text{otherwise}
\end{cases}}\\
&={\begin{cases} \ch_{Z(\alpha,s(\alpha))}([\alpha, s(\alpha),\res_{\alpha}\circ\ind_{\xi}\circ \res_{\zeta} x][\xi,\zeta,x]) &\text{if }\ind_{\xi}\circ \res_{\zeta} x\in Z(s(\alpha))\\
0 &\text{otherwise}
\end{cases}}\\
\intertext{For $[\alpha, s(\alpha),\res_{\alpha}\circ\ind_{\xi}\circ \res_{\zeta} x],[\xi,\zeta,x]$ to be composable $r(\xi)=s(\alpha)$ and so}
&={\begin{cases} \ch_{Z(\alpha,s(\alpha))}([\alpha\xi,\zeta,x]) &\text{if }\ind_{\xi}\circ \res_{\zeta} x\in Z(s(\alpha))\\
0 &\text{otherwise}
\end{cases}}\\
&={\begin{cases} 1 &\text{if }[\alpha\xi,\zeta,x]=[\alpha,s(\alpha),x]\\
0 &\text{otherwise}.
\end{cases}}
\end{align*}
Now $[\alpha\xi,\zeta,x]=[\alpha,s(\alpha),x]$ implies $\ind_{\alpha\xi}\circ \res_{\zeta}=\ind_{\alpha}$ on some neighborhood, so 
$\ind_{\xi}\circ\res_\zeta=\Id_{Z(s(\alpha))}$ on that same neighborhood.  There exists $a,b\in \mathcal{B}$ such that $F(\alpha\xi)a=F(\alpha)b$ and $F(\zeta)a=b$ so $F(\alpha\xi)a=F(\alpha\zeta)a$.  Thus, since morphisms in  $\mathcal{B}$ are monic, we have $F(\xi)a=F(\zeta)a$.  That is $[\xi,\zeta,x]=[s(\alpha),s(\alpha), x]$.  So that $\ch_{Z(s(\alpha),\alpha)}*\ch_{Z(\alpha,s(\alpha))}=\ch_{Z(s(\alpha),s(\alpha))}$.
$\Box$ \medskip

Before proving that $\Upsilon$ is an isomorphism we need a lemma. 

\begin{lem}For $\alpha,\beta\in\E$
\label{lem:Z int}
\begin{enumerate}
\item \label{it: 1}if $r(\alpha)=r(\beta)=X$ and $Z(\alpha)=Z(\beta)$ then $S_\alpha S_\alpha^*=S_\beta S_\beta^*$;
\item \label{it: 2} if $Z(\alpha)\cap Z(\beta)=\emptyset$ then $S_\alpha^*S_\beta=0$.
\end{enumerate}
\end{lem}

\noindent {\sc Proof.} For Item \eqref{it: 1}, note that $F(\alpha),F(\beta)$ is a cospan in $\mathcal{B}$ so by right Ore there exist $a,b\in \mathcal{B}$ such that $F(\alpha)a=F(\beta)b$.  If $\eta\in \E$ such that $r(\eta)=s(\alpha)$ and $F(\eta)=a$ then there exists an infinite path $x\in Z(\alpha)=Z(\beta)$ such that $x(F(\alpha)a)=\alpha\eta$. But then   $\alpha\eta=x(F(\alpha)a)=x(F(\beta)b)=\beta x_2(F(\beta),b)$; thus $\{\alpha\eta: F(\eta)=a, r(\eta)=s(\alpha)\}=\{\beta\gamma: F(\gamma)=b, r(\gamma)=s(\beta)\}$.  Therefore
\[
s_\alpha s_\alpha^*=\sum_{F(\eta)=a, r(\eta)=s(\alpha)}s_{\alpha\eta}s_{\alpha\eta}^*=\sum_{F(\gamma)=b, r(\gamma)=s(\beta)}s_{\beta\gamma}s_{\beta\gamma}^*=s_\beta s_\beta^*.
\]

For Item \eqref{it: 2}, again $F(\alpha),F(\beta)$ is a cospan in $\mathcal{B}$ so by right Ore there exist $a,b\in \mathcal{B}$ such that $F(\alpha)a=F(\beta)b$.  So 
\[
s_\alpha^*s_\beta=\sum_{F(\eta)=a, r(\eta)=s(\alpha), F(\gamma)=b, r(\gamma)=s(\beta)} s_\eta s_{\alpha\eta}^*s_{\beta\gamma}s_\gamma^*.
\]
It suffices to see $\{\alpha\eta: F(\eta)=a, r(\eta)=s(\alpha)\}\cap \{\beta\gamma: F(\gamma)=b, r(\gamma)=s(\beta)\}=\emptyset$.  But if there exist $\alpha\eta=\beta\gamma$ then for any $y\in Z(s(\eta))$, $\ind_{\alpha\eta} y=\ind_{\beta\gamma} y\in Z(\alpha)\cap Z(\beta)$ a contradiction. $\Box$\medskip

The proof of the following is modified from Theorem~4.2 in \cite{KPRR}.

\begin{thm} Let $F:\E\to \B$ be a row-finite Kumjain Pask fibration. If all morphisms of $\cal B$ are monic, then
$\Upsilon:C^*(F)\to C^*(G_F)$ characterized by $s_\alpha\mapsto \ch_{Z(\alpha,s(\alpha))}$ is an isomorphism.
\end{thm}

\noindent {\sc Proof.}
To show $\Upsilon$ is an isomorphism we construct an inverse for it.  We do this locally, starting with $C(Z(X))$ (for a each object $X$) and then gluing the resulting maps together.  Notice that the map $\phi_{\alpha,\beta}:x\mapsto [\alpha,\beta,\ind_\beta x]$ is a homeomorphism of $Z(s(\beta))$ to $Z(\alpha,\beta)$. Also for an object $X$ in $\E$, $C(Z(X))=\cspan\{\ch_{Z(\alpha)}:r(\alpha)=X\}$ by the Stone-Weierstrass Theorem.

By Lemma~\ref{lem:Z int} Item~\eqref{it: 1}, $\theta_X': \{\ch_{Z(\alpha)}:r(\alpha)=X\}\to C^*(F)$ given by $\theta_X'(\ch_{Z(\alpha)})=S_\alpha S_\alpha^*$ is well-defined.  Extend $\theta_X'$ linearly to the span of the $\ch_{Z(\alpha)}$'s.  We want to extend $\theta_X'$ to $C(Z(X))$ but to do this we need to see $\theta_X'$ is norm decreasing.  

Since each $Z(\alpha)$ is compact open and $F$ is row-finite we can use Lemma~\ref{lem:Z int} Item~\eqref{it: 2} to see any element in $\spn\{\ch_{Z(\alpha)}:r(\alpha)=v\}$ can be written as $\sum_{\alpha\in I} r_\alpha \ch_{Z(\alpha)}$ such that for $\alpha\neq \beta\in I$, $Z(\alpha)\cap Z(\beta)=\emptyset.$  Thus $\theta_X'(\sum_{\alpha\in I} r_\alpha \ch_{Z(\alpha)})=\sum_{\alpha\in I} r_\alpha s_\alpha s_\alpha^*$ where $\{s_\alpha s_\alpha^*\}_{\alpha\in I}$ is a set of mutually orthogonal projections.  Therefore $\|\theta_X'(\sum_{\alpha\in I} r_\alpha \ch_{Z(\alpha)})\|=\max\{|r_\alpha|\}=\|\sum_{\alpha\in I} r_\alpha \ch_{Z(\alpha)}\|$ and so $\theta_X'$ is norm preserving and thus extends to a linear map $\theta_X: C(Z(X))\to C^*(F)$; in fact it is a $*$-homomorphism.

Define  $\rho_{\alpha,\beta}: C(Z(\alpha,\beta))\to C^*(F)$ by $f\mapsto \theta_{s(\alpha)}\circ \phi_{\alpha,\beta}\inv(f)$. So that 
\[
\rho_{\alpha,\beta}(\ch_{Z(\alpha\gamma,\beta\gamma)})=\theta_{s(\alpha)}(Z(\gamma))=s_\gamma s_\gamma^*.
\]
We claim that 
\begin{equation}
\label{eq1}
\rho_{\alpha,\beta}(f)=\sum_{r(\gamma)=s(\alpha), F(\gamma)=c} s_\gamma \rho_{\alpha\gamma,\beta\gamma}(f|_{Z(\alpha\gamma,\beta\gamma)})s_\gamma^*.
\end{equation}
Since $\rho_{\alpha\gamma,\beta\gamma}$, $\rho_{\alpha,\beta}$ and $\text{Ad}_{s_{\gamma}}$ are linear and continuous for all $\gamma$ and $C(Z(\alpha,\beta))=\cspan\{\ch_{Z(\alpha\eta,\beta\eta)}\}$, it suffices to check \eqref{eq1} on characteristic functions of the form $\ch_{Z(\alpha\eta,\beta\eta)}$. Choose such an $\eta$. By the right Ore there exists $a,b$ such that $F(\beta)F(\eta)a=F(\beta)cb.$  Since  morphisms are monic we get $F(\eta)a=cb$ so that $F(\alpha)F(\eta)a=F(\alpha)cb$ as well.  Therefore 
\[
\ch_{Z(\alpha\eta,\beta\eta)}=\sum_{F(\xi)=b, F(\gamma)=c, \eta\zeta=\gamma\xi} \ch_{Z(\alpha\gamma\xi,\beta\gamma\xi)}=\sum_{F(\gamma)=c} \sum_{F(\xi)=b,\eta\zeta=\gamma\xi} \ch_{Z(\alpha\gamma\xi,\beta\gamma\xi)}.
\]  
Notice that $\ch_{Z(\alpha\eta,\beta\eta)}|_{Z(\alpha\gamma,\beta\gamma)}=\sum_{F(\xi)=b,\eta\zeta=\gamma\xi} \ch_{Z(\alpha\gamma\xi,\beta\gamma\xi)}$.
Now 
\begin{align*}
\rho_{\alpha,\beta}(\ch_{Z(\alpha\eta,\beta\eta)})&=\rho_{\alpha,\beta}\biggl(\sum_{F(\gamma)=c} \biggl(\sum_{F(\xi)=b,\eta\zeta=\gamma\xi} \ch_{Z(\alpha\gamma\xi,\beta\gamma\xi)}\biggr)\biggr)\\
&=\sum_{F(\gamma)=c} \biggl(\sum_{F(\xi)=b,\eta\zeta=\gamma\xi} s_{\gamma\xi}s_{\gamma\xi}^*\biggr)\\
&=\sum_{F(\gamma)=c}s_\gamma\biggl( \sum_{F(\xi)=b,\eta\zeta=\gamma\xi} s_{\xi}s_{\xi}^*\biggr)s_\gamma^*\\
&=\sum_{F(\gamma)=c}s_\gamma( \sum_{F(\xi)=b,\eta\zeta=\gamma\xi} \rho_{\alpha\gamma,\beta\gamma}(\ch_{Z(\alpha\gamma\xi,\beta\gamma\xi)}))s_\gamma^*\\
&=\sum_{F(\gamma)=c}s_\gamma(\rho_{\alpha\gamma,\beta\gamma}(\ch_{Z(\alpha\eta,\beta\eta)}|_{Z(\alpha\gamma,\beta\gamma)}))s_\gamma^*
\end{align*}
as desired.

For $f\in C_c(G)$ using Proposition~\ref{disjoint basis} we can write $\supp(f)$ as a disjoint union of $Z(\alpha^i, \beta^i)$.  We define
\[
\theta(f)=\sum s_{\alpha^i}(\rho_{\alpha^i,\beta^i} (f|_{Z(\alpha^i,\beta^i)})))s_{\beta^i}^*=\sum s_{\alpha^i}(\theta_{s(\alpha_i)}(\phi_{\alpha^i,\beta^i}\inv (f|_{Z(\alpha^i,\beta^i)})))s_{\beta^i}^*.
\]
We need to show $\theta(f)$ is well-defined; that is does not depend on the decomposition of the support of $f$.  Suppose $\supp(f)=\bigcup Z(\sigma^j,\tau^j)$.  Then $\supp(f)=\cup (Z(\alpha^i,\beta^i)\cap Z(\sigma^j,\tau^j))$.  Since each of these intersections is a union of sets of the form $Z(\alpha^i\gamma,\beta^i\gamma)$ with $F(\gamma)=a_i$, it suffices to show that for $f$ with $\supp(f)\subset Z(\alpha,\beta)$ that
\begin{align*}
s_\alpha (\rho_{\alpha,\beta}(f)) s_\beta^*&=\sum_{r(\gamma)=s(\alpha), F(\gamma)=b} s_{\alpha\gamma}(\rho_{\alpha\gamma,\beta\gamma}(f|_{Z(\alpha\gamma,\beta\gamma)})))s_{\beta\gamma}^*\\
&=s_\alpha\bigl(\sum_{r(\gamma)=s(\alpha), F(\gamma)=b} s_{\gamma}(\theta_{s(\gamma)}(\phi_{\alpha\gamma,\beta\gamma}\inv(f|_{Z(\alpha\gamma,\beta\gamma)})))s_{\gamma}^*\bigr)s_\beta^*\\
\end{align*}
But we know that $\rho_{\alpha,\beta}(f))=\sum_{r(\gamma)=s(\alpha), F(\gamma)=b} s_{\gamma}(\theta_{s(\gamma)}(\phi_{\alpha\gamma,\beta\gamma}\inv(f|_{Z(\alpha\gamma,\beta\gamma)})))s_{\gamma}^*$ from above.  Thus $\theta$ is well-defined.

To show that $\theta: C^*(G_F)\to C^*(F)$ is continuous it is enough to show $\theta: C_c(G_F)\to C^*(F)$ is continuous in the inductive limit topology.  Suppose $f_i\to f$ uniformly with $\supp(f_i)$ eventually contained in a fixed compact set $K$.  Let $K\subset \cup Z(\alpha^j, \beta^j)$. Since $\supp(f_i)\subset K\subset \cup Z(\alpha^j, \beta^j)$, there exists $N$ such that for $i\geq N$ such that $\theta(f_i)=\sum s_{\alpha^j}\rho_{\alpha^j,\beta^j} (f_i|_{Z(\alpha^j,\beta^j)})s_{\beta^j}^*$.  But the maps $\rho_{\alpha^j,\beta^j}$ are all norm decreasing for the uniform norm on $C(Z(\alpha^i,\beta^i))$ and the $C^*$-norm on $C^*(F)$.  Thus $\theta(f_i)\to \theta(f)$ and hence $\theta$ extends to a continuous map $\theta: C^*(G_F)\to C^*(F)$.

It remains to show that $\theta$ is a $*$-homomorphism.  Since $\theta$ is continuous and linear and $C^*(G_F)$ is the closed span of $\{\ch_{Z(\alpha,\beta)}\}$  it suffices to show that $\theta$ is multiplicative on these functions; but

\begin{align*}
\theta(\ch_{Z(\alpha,\beta)}*\ch_{Z(\mu,\nu)})&=\theta(\ch_{Z(\alpha,\beta)Z(\mu,\nu)})\\
&=\theta\left(\sum_{\beta\eta=\mu\lambda, F(\eta)=a, F(\lambda)=b} \ch_{Z(\alpha\eta,\nu\lambda)}\right)\\
&=\sum_{\beta\eta=\mu\lambda, F(\eta)=a, F(\lambda)=b} \theta(\ch_{Z(\alpha\eta,\nu\lambda)})\\
&=\sum_{\beta\eta=\mu\lambda, F(\eta)=a, F(\lambda)=b} s_{\alpha\eta} s_{\nu\lambda}^*\\
&=s_\alpha s_\beta^* s_\mu s_\nu^*=\theta(\ch_{Z(\alpha,\beta)})\theta(\ch_{Z(\mu,\nu)}).\end{align*}

Also,

\[\theta(\ch_{Z(\alpha,\beta)})^* = (s_\alpha s_\beta^*)^* = s_\beta s_\alpha^*=\theta(\ch_{Z(\alpha,\beta)}^*) .
\]

Lastly, by the definition of $\theta$, it is an inverse for $\Upsilon$.
$\Box$ \medskip

It is possible that $Z(\alpha,s(\alpha))=Z(\beta,s(\beta))$ with $\alpha\neq \beta$.  But for this to be true there must exist an $a$ such that $F(\alpha)a=F(\beta)a$.   For $[\alpha, s(\alpha), x]\in Z(\beta,s(\beta))$ we must have $[\alpha,s(\alpha), x]=[\beta, s(\beta),x]$ so that for $\gamma=x(a)$,$\alpha\gamma=\ind_{\alpha}(x)(F(\alpha)a)=\ind_{\beta}(x)(F(\beta)a)=\beta\gamma$.    Thus
\[
s_\alpha=s_\alpha\sum_{F(\gamma)=a, r(\gamma)=s(\alpha)}s_\gamma s_\gamma^*= \sum_{F(\gamma)=a, r(\gamma)=s(\alpha)}s_{\alpha\gamma} s_\gamma^*=\sum_{F(\gamma)=a, r(\gamma)=s(\beta)}s_{\beta\gamma} s_\gamma^*=s_\beta.
\]
Now if $\mathcal{B}$ is right cancellative  $Z(\alpha,s(\alpha))=Z(\beta,s(\beta))$ implies $F(\alpha)a=F(\beta)a$ so $F(\alpha)=F(\beta)$ and this gives $\alpha=\beta$.

\begin{proposition}
\label{prop:inj}
Let $F:\mathcal{E}\to \mathcal{B}$ be a KPf with $\mathcal{B}$ left cancellative.  The map $\alpha\mapsto s_\alpha$ from $\mathcal{E}\to C^*(F)$ is injective if and only if $\mathcal{B}$ is right cancellative.
\end{proposition}

\noindent{\sc Proof.}
Let $\Upsilon:C^*(F)\to C^*(G_F)$ be the $*$-homomorphism from Lemma~\ref{hom} which sends $s_\alpha\mapsto \ch_{Z(\alpha,s(\alpha))}$.  Suppose $s_\alpha=s_\beta$, then $\ch_{Z(\alpha,s(\alpha))}=\Upsilon(s_\alpha)=\Upsilon(s_\beta)=\ch_{Z(\beta,s(\beta))}$.  Thus 
\[Z(\alpha,s(\alpha))=Z(\beta,s(\beta)).\]
  Chose $x\in Z(s(\alpha))$.  Then $[\alpha, s(\alpha),x]=[\beta,s(\beta),x]$.  In particular
  
\begin{align*}
\ind_\alpha(x)&=\ind_\beta(x), \text{~and} \\
\intertext{there exist $a,b\in \mathcal{B}$ with}
F(\alpha)a&=F(\beta)b, \text{~and}\\
F(s(\alpha))a&=F(s(\beta))b\\
\intertext{which implies}
F(\alpha)a&=F(\beta)a.
\end{align*}
  
  If $\mathcal{B}$ is right cancellative, we have $F(\alpha)=F(\beta)$ and so $\alpha=\ind_\alpha(x)(F(\alpha))=\ind_\beta(x)(F(\beta))=\beta$.  That is the map $\alpha\mapsto s_\alpha$ is injective.
  
Now suppose $\alpha\neq \beta$ and $s_\alpha=s_\beta$  then $F(\alpha)\neq F(\beta)$: otherwise $0\neq s_\alpha^*s_\alpha=s_\alpha^*s_\beta=0$.  By above we have there exists $a\in \mathcal{B}$ with $F(\alpha)a=F(\beta)a$ so $\mathcal{B}$ is not right cancellative. $\Box$
\medskip

Since every $C^*$-algebra can be represented faithfully on a Hilbert space the previous result says row-finite Kumjian-Pask fibration $F:\E\to \B$ with $\B$ left cancellative can be represented faithfully on a Hilbert space if and only if the base category is also right cancellative.  In special circumstances the infinite path representation gives a faithful representation of the Kumjian-Pask fibration.

\begin{defn}
We say $x\in F^\infty$ is aperiodic if for all $\alpha,\beta\in \mathcal{E}$, $\res_\alpha(x)=\res_\beta(x)$ imples $\alpha=\beta$.  We say $F:{\cal E}\rightarrow {\cal B}$ is aperiodic if $Z(X)$ contains an aperiodic path for every $X\in \text{Obj}(\mathcal{E})$.
\end{defn}

For a groupoid $G$ and $x\in G^{(0)}$ denote $xGx=\{\gamma\in G:r(\gamma)=s(\gamma)=x\}$.

\begin{thm} \label{aper_implies_top_prin}
If $\mathcal{E}$ is aperiodic then $G_F$ is topologically principal.
\end{thm}

\noindent{\sc Proof.}
We need to show $T=\{x\in F^\infty: xG_Fx=\{x\}\}$ is dense in $F^\infty$.  Let $Z(\alpha)$ be a basic open set in $F^\infty$, it suffices to show $T\cap Z(\alpha)\neq \emptyset$.  Now $Z(\alpha)=Z(\alpha,s(\alpha))\cdot Z(s(\alpha))$ so it suffices to show $Z(s(\alpha))\cap T=\emptyset$, since if $x\in Z(s(\alpha))\cap T$ then $[\alpha,s(\alpha),x]\cdot x\in Z(\alpha)\cap T$.
By assumption there exists an aperiodic path $x\in Z(s(\alpha))$.  We need to show $xG_F x=\{x\}$.  Now if $[\alpha,\beta,x]\in xG_Ex$ then $\ind_\alpha\circ \res_\beta(x)=x$ that is $\res_\beta(x)=\res_\alpha(x)$. Since $x$ aperiodic we have $\alpha=\beta$ so that $[\alpha,\beta,x]=[r(x),r(x),x]$.  That is $x\in T\cap Z(s(\alpha))$ and so $G_F$ is topologically principal as desired.
$\Box$
\medskip

Let $I_\lambda$ be the ideal such that $C^*(G_F)/I_\lambda=C^*_r(G_F)$.  If  $G_F$ is amenable then $C^*(G_F)\cong C^*_r(G_F)$. (Amenability is extremely complicated for groupoids.  For a full discussion see \cite{ADR}. Fortunately here we only need that $G_F$ amenable implies $C^*(G_F)\cong C^*_r(G_F)$.)

\begin{cor}
\label{CK}
Let $F:\mathcal{E}\to \mathcal{B}$ be an aperiodic Kumjian-Pask fibration over a right cancellative category $\cal B$.  
\begin{enumerate}
\item Then the kernel of the infinite path representation of $C^*(G_F)$ is contained in $I_\lambda$.  In particular $\mathcal{E}$ represents faithfully on $\ell^2(F^\infty)$.
\item If additionally $G_F$ is amenable and $\pi: C^*(F)\to B$ is a $*$-homomorphism.  Then $\pi$ is faithful if and only if $\pi(P_X)\neq 0$ for all objects $X$ in $\E$.
\end{enumerate}
\end{cor}

\noindent{\sc Proof.}
By  Theorem~5.5 in \cite{BCFS} if $G_F$ is topologically principal  then every nonzero ideal of $C_r^*(G_F)$ intersects $C_0(G_F^0)$ nontrivially.
 
For Item~(1), let $q: C^*(G_F)\to C^*_r(G_F)$ be the quotient map and $J$ be the kernel of the infinite path representation.  Then $q(J)$ is an ideal in $C^*_r(G_F)$.  Now $C_0(G_F^{(0)})\cap J=\{0\}$ and since $q|_{C_0(G_F^{(0)})}$ is the identity we have $q(J)\cap C_0(G_F^{(0)})=\{0\}$ which implies $q(J)=\{0\}$ hence $J\subset I_\lambda$ as desired.  For the second statement $\Upsilon(S_\alpha)=\ch_{Z(\alpha,s(\alpha))}\in C_c(G_F)$ and $C_c(G_F)\cap I_\lambda=\{0\}$.

For Item~(2), if $G_F$ amenable $C^*(G_F)\cong C^*_r(G_F)$ so any ideal of $C^*(G_F)\cong C^*(F)$ has nontrivial intersection with $C_0(G_F^{(0)})$.  In particular $\ker(\pi)$ contains a $\ch_{Z(\alpha)}=\ch_{Z(\alpha,s(\alpha))}*\ch_{Z(\alpha,s(\alpha))}^*=\Upsilon(S_\alpha S_\alpha^*)$ for some $\alpha\in \E$.  Thus $\pi(S_\alpha)=0$ which gives $\pi(P_{s(\alpha)})=\pi(S_\alpha^*S_\alpha)=0$.
$\Box$
\medskip

Corollary~\ref{CK} Item~(2) is typically referred to as the Cuntz-Krieger uniqueness theorem.  Observe that it depends on the amenability of $G_F$.    For a $k$-graph, $F=d: \Lambda\to \N^k$, the amenability of $G_F$ was shown in Theorem~5.5 in \cite{KP}.  On the other hand, if the Kumjian-Pask fibration is $F=\Id_H: H\to H$ for $H$ a nonamenable group (say ${\mathbb{F}}_2)$, then $G_F=H$. Thus we can not expect $G_F$ to be amenable in general.  It would thus be desirable to find a sufficient condition, or ideally necessary and sufficient condition, on $F$ which would ensure the amenability of $G_F$.
 
Now suppose $F: \E\to \B$ with $\B$ an abelian monoid.  Then $F$ induces a map $c: G_F\to H$ where $H$ is the Grothendieck group constructed from $\B$.   In this case Theorem~9.3 of \cite{S} shows that $G_F$ is amenable if  $c\inv(0)$ is an amenable groupoid.  In \cite[Theorem~9.8]{S}, Spielberg shows that under certain conditions for a category of paths, $c\inv(0)$ is an AF groupoid and hence amenable.  Although the presence of multiple objects in $\B$ will complicate the matter, we expect a similar result to hold in the general case.

Besides the investigation of amenability for general Kumjian-Pask fibrations, and consideration of the properties of specific examples and classes of examples of $C^*$-algebras arising from the constructions in the present paper, two other avenues of investigation present themselves -- consideration of the construction in the absence of either the countability condition on the slice categories $\B / B$ or row-finiteness.

In the latter case, it would appear the correct generalization of Definition \ref{CKsystem} Item (6) would be:  For all $X$, and for all morphisms $b:B\rightarrow F(X)$ in ${\cal B}$, and all  $\alpha \in X\E \cap F^{-1}(b)$

\[  S_\alpha S^*_\alpha \leq P_X.\]

\end{document}